\documentclass[11pt]{article}
\usepackage{amsfonts}
\usepackage{mathrsfs}
\usepackage{amsthm}
 \usepackage{graphicx} %We can use any anther package if necessary
\usepackage{amsmath} % If you want to use AMS-LaTeX, comment out there
\usepackage{amssymb} % Two commands
\usepackage{enumerate}
\usepackage{epsfig}
\usepackage{color}

\newtheorem{thm}{\hskip\parindent Theorem}
\newtheorem{lem}{\hskip\parindent Lemma}

\numberwithin{equation}{section}

\def\Ai{\mathrm {Ai}}
\def\Bi{\mathrm {Bi}}
\def\Re {\mathop{\rm Re}\nolimits}
\def\Im {\mathop{\rm Im}\nolimits}
\def\Res {\mathop{\rm Res}\limits}

\title{Weights with both absolutely  continuous and discrete components:  Asymptotics via  the Riemann-Hilbert approach  }
\author{Xiao-Bo Wu$^{\,1}$, Yu Lin$^{\,2}$, Shuai-Xia Xu$^{\,3}$ and  Yu-Qiu Zhao$^{\,1}$\footnote{Corresponding author (Yu-Qiu Zhao).
 {\it{E-mail
address:}} {stszyq@mail.sysu.edu.cn} } }

\date{
\hbox{\small \emph{$^1$ Department of Mathematics, Sun Yat-sen University, Guangzhou,
China}}
 \hbox{\small
\emph{$^2$ Department of Mathematics, South China University of Technology, Guangzhou, China}
 }
 \hbox{\small
\emph{$^3$ Institut Franco-Chinois de I'Energie Nucl\'eaire, Sun Yat-sen University, Guangzhou,
China}}}

\begin{document}
\maketitle

\begin{abstract}
{We study the uniform asymptotics for the orthogonal polynomials with respect to weights composed of both absolutely continuous measure and discrete measure, by taking a special class of the sieved Pollazek Polynomials as an example.  The Plancherel-Rotach type asymptotics  of the   sieved Pollazek Polynomials are obtained  in the whole complex plane.  The Riemann-Hilbert method is applied to derive the results.  A main  feature of the treatment  is the appearance of a new band
consisting  of two adjacent intervals, one of which is  a portion of  the support of the absolutely continuous measure,   the other is the discrete band.
}

\vspace{.4cm}

\textbf{Keywords:}\;Uniform asymptotics; discrete orthogonal  polynomials; Riemann-Hilbert approach;  sieved Pollaczek polynomials; Airy function.

\textbf{Mathematics Subject Classification 2010}: 41A60, 33C10, 33C45

\end{abstract}

%====================================================================================%

\section{Introduction}  \label{introduction}
The method of Deift and Zhou has found various applications in the asymptotic studies of orthogonal polynomials.
The first few examples, published in 1999, are polynomials  with absolutely  continuous weights; see, e.g.,  Deift {\it et al.} \cite{dkmv2}. The   powerful method is based on the Riemann-Hilbert problem (RHP) formulation  of the orthogonal polynomials observed by Fokas, Its and Kitaev \cite{Fokas}. A crucial idea is a deformation of the contours associated with  the factorization of the oscillating jump matrices.
Technique difficulties usually lie in the  construction of the parametrices at critical points.

In 2007, Baik {\it et al.}\;\cite{bkmm}  developed a general method for  the asymptotics of discrete orthogonal polynomials  by using the Riemann-Hilbert approach. The starting point of their  investigation  is  the  interpolation problem (IP)  for    discrete orthogonal polynomials, introduced in  Borodin and Boyarchenko   \cite{Borodin}.
 A key step is to turn the IP  into a RHP, and then the  Deift-Zhou method for oscillating RHP may  apply.
For the case when all nodes are real,  the real line is divided into intervals termed void, saturated region, or band, associated with the equilibrium measure; see \cite{bkmm,Bleher} when the nodes are regularly distributed.
In general case, an open interval is called a saturated region, if the ratio of the density of the polynomial zeros and the density of the nodes is $1$, called void if the ration is $0$, and band otherwise.

Since then, much attention has been attracted to this topic.   For example, in an attempt to achieve global asymptotics, with global referring to the   domains of uniformity,
   Wong  and coworkers considered  cases with finite nodes \cite{DaiWong,LinWong,LinWongC},  and infinite nodes \cite{OuWong,WangWong}  regularly distributed. Very recently, the   present authors \cite{WuLinXuZhao} have studied  the uniform asymptotics for discrete orthogonal  polynomials on infinite nodes with an accumulation point, the mass showing  a singular behavior there.

A major modification to the method has been  made by
Bleher and Liechty \cite{BleherE,Bleher}   in the treatment of the    band-saturated region endpoints.  The example they take  is  a system of discrete orthogonal polynomials
with respect to a varying exponential weight   on a regular infinite lattice.  Here  regular infinite lattice means that the infinite nodes are equally spaced.

Other than the discrete orthogonal polynomials and those with absolutely continuous weights, there are interesting  mixed-type cases that the orthogonal measures are supported on both intervals and discrete nodes. Examples can be found in Askey and Ismail \cite[Ch. 6]{AskeyIsmail}, and in Ismail \cite[p. 156]{Ismail}, where the random walk polynomials are shown having this feature.   A natural question  arises here: How could the Riemann-Hilbert approach be used to handle such problems?

In this paper, we illustrate the method by taking as an example a class of sieved Pollaczek polynomials.
  A significant fact is  that the corresponding orthogonal measure  consists of an absolutely continuous part on $[-1, 1]$, and a discrete part having infinite many mass points with the endpoint $1$ as an accumulation point;
see Charris and ismail \cite{charris-ismail}  and Wang and Zhao \cite{wang-zhao}.

It is known that the sieved Pollaczek polynomials  $p_n(z)$   possess the orthogonal relation
\begin{equation}\label{orthogonal-relation}
\int_{-1}^{1} p_n(x)p_m(x)w_c(x)dx+\sum_{k=1}^{\infty}p_n(x_k)p_m(x_k)w_d(x_k)=\delta_{n,m},~~n,m=0,1,2,\cdots ,
\end{equation}
 where the absolutely  continuous weight
\begin{equation}\label{weight-continuos}
    w_c(x)=\frac{2\sin \theta}{\pi}\exp\left \{\frac{(\pi-2\theta)b}{\sin \theta}\right \} \; \left |\Gamma\left (1+\frac{ib}{\sin \theta}\right )\right |^2=\frac{2b\exp\left \{ \frac {(\pi-2\theta) b}{\sin\theta}\right \}} {\sinh\frac {b\pi}{\sin\theta}}
\end{equation}for $\theta=\arccos x\in(0,\pi)$ while $x\in (-1, 1)$, and the
  mass
\begin{equation}\label{discontinuous-weight}
  w_d(x_k)=\frac{4b^3}{k^3\sqrt{1+b^2/k^2}}(\sqrt{1+b^2/k^2}-b/k)^{2k}
\end{equation}
with nodes $x_k=\sqrt{1+b^2/k^2}$, $k=1,2,\cdots$; see \cite{charris-ismail,Ismail, wang-zhao} for detailed determination of the orthogonal measure.   It is easily seen that $w_d(x_k)\sim 4b^3 e^{-2b}/k^3$ as $k\to\infty$, which indicates a singularity of the discrete weight at the accumulation point $x=1$.

The polynomials $p_n(z)=p_n(z;b)$ can also be defined by the three-term recurrence relation
\begin{equation}\label{pn}
p_{n+1}(z)+p_{n-1}(z)+\frac{2b}{n+1}p_n(z)=2zp_n(z), \quad n=1,2,\cdots
\end{equation}
with initial values $p_0(z)=1$ and $p_1(z)=2z-2b$; see Charris and ismail \cite{charris-ismail}, Ismail\cite{Ismail}, and Wang and Zhao \cite{wang-zhao} and the references therein. It is seen from \eqref{pn} and the corresponding initial conditions that
\begin{equation}\label{-b-b}
p_n(z;-b)=(-1)^np_n(-z;b).
\end{equation}Hence, without loss of generality, we may assume that $b>0$.

Methods other than the Riemann-Hilbert approach may be applied to obtain asymptotics of the  sieved Pollaczek  polynomials; cf., e.g., Szeg\"{o} \cite{Szego1975} and Wong and Zhao \cite{WongZhao}.  Indeed, in an earlier work, Wang and Zhao \cite{wang-zhao} have considered the uniform asymptotic expansions for the  polynomials on the real line, in particular at the turning points $x=-1+b/n$ and $x=1+b/n$ and the endpoints $\pm 1$,
 by using an integral method.  The expansion  in an $O(1/n)$ neighborhood of $-1+b/n$ is in terms of the Airy function, while at   $1+b/n$, where the polynomials oscillate on both sides, we need a combination of the Airy functions to describe the behavior.    The asymptotics of the extreme    zeros are also obtained in \cite{wang-zhao}. However, the derivation is limited to the real line, and, rigorously, there are gaps between the intervals of  uniformity. Further study is desirable for such polynomials.

The main purpose of the present investigation is to derive uniform  asymptotic approximations on the whole complex  plane   for the orthogonal  polynomials with
  both absolutely continuous measure and discrete measure, using   the Riemann-Hilbert approach and taking the  sieved Pollaczek polynomials as an example.

To this aim,
   first we formulate the mixed Riemann-Hilbert and interpolation problem for the polynomials. Then, we convert the problem into a Riemann-Hilbert problem by using  the notion of band and saturated region of  Baik {\it et al.} \cite{bkmm}, and the treatment of the band-saturated region endpoints by  Bleher and Liechty \cite{BleherE, Bleher}. The Deift-Zhou nonlinear steepest descent method for oscillating Riemann-Hilbert problems    plays a central part from then on.  The main idea here is the oscillating contour consists of the  interval of  the absolutely continuous measure,  joined by the discrete band.

 Technically, there are several facts in the  analysis worth mentioning. Several   auxiliary functions    $d_E$,  $d_I$ and $\chi$
are introduced at an early stage $Y\to U$, to simplify the jump conditions, and to clarify  the construction of the outer parametrix for $N$.  The $g$-function is supported on an infinite interval, so that the contours are finite later in the RHP for $T$.   In the    parametrix for $N$, we bring in extra singularities to fit the matching conditions on the boundaries of  the shrinking neighborhoods of the MRS numbers $\alpha$ and $\beta$, in  which the local papametries are constructed. The  phase condition  on the band $(\alpha, \beta)$ plays a role in the determination of the equilibrium measure, and  very careful estimates of the $\phi$-functions are also needed since the  domains of   local parametries are   shrinking.

The paper is arranged as follows. In Section \ref{statement-of-result}, we state the main asymptotic approximations in regions covering the upper half plane.  In Section \ref{RHP-IP}, as the starting point of our analysis,
we formulate  the Riemann-Hilbert and interpolation    problem (RHP and IP)  for the sieved Pollaczek polynomials. In Section \ref{removing-poles}, by  removing  the poles at the nodes, the problem is  turned into a RHP  for a matrix function $U(z)$. In Section \ref{MRS}, we calculate the MRS numbers and bring in auxiliary functions such as the $g$-function and $\phi$-functions.   The nonlinear steepest descent analysis is carried in  Section \ref{steepest-descent-analysis}. The proof of the main asymptotic theorem is provided in  Section \ref{proof-of-theorem} by using the analysis in previous sections. Several asymptotic quantities are calculated in Section \ref{comparison}, and a comparison of results is made with   the known ones in \cite{wang-zhao}.

\section{Main results:  Uniform  asymptotic approximations }
 \setcounter{equation} {0}  \label{statement-of-result}

We derive asymptotic approximations for the orthonormal
sieved Pollaczek
 polynomials $p_n(z)$ in overlapped
domains covering the whole complex plane, based on the Riemann-Hilbert analysis
carried out. In view of the symmetry with respect to the real line, we need only to work on
the upper half-plane.

To describe the results, we   introduce several constants and auxiliary  functions.

The soft edges are  located at the MRS numbers  $\alpha=\alpha_n$ and $\beta=\beta_n$. A detailed analysis  of these numbers   will be carried out  in   Section \ref{MRS}. We will see that  $
 \alpha =-1+\frac b n+O\left(\frac{1}{n^2}\right)$ and  $\beta  = 1+\frac b n+O\left(\frac{1}{n^2}\right)$ as $n\to\infty$.

Next, we define the functions
\begin{equation*} \gamma(z)=\frac{b}{\sqrt{z^2-1}}
 \quad \mbox{and}\quad \varrho(z)=\left ( \frac {z-\beta}{z-\alpha} \right )^{1/4},
\end{equation*}analytic respectively in $\mathbb{C}\setminus [-1, 1]$ and $\mathbb{C}\setminus [\alpha, \beta]$, where the branches are chosen such that $\arg (z-\alpha), ~\arg(z-\beta), ~\arg (z\pm 1)\in (-\pi, \pi)$; cf. Section \ref{removing-poles} and Section \ref{steepest-descent-analysis} below.
We also need  the scalar function
\begin{equation}\label{phi-0}
 \phi_0(z) =\frac {e^{\pm \pi i}} {\varphi(z)} \frac {\Gamma(1-\gamma(z))}{\Gamma(1+\gamma(z))} e^{-2\gamma(z)+2\gamma(z) \ln (-\gamma(z))\mp \pi i (\gamma(z) +1/2)}~~\mbox{for}~~\pm \Im z>0,
\end{equation}
where $\varphi(z)=z+\sqrt{z^2-1}$ is analytic in $\mathbb{C}\setminus [-1, 1]$ such that $\varphi(z)  \approx 2z$ for large $z$, and the logarithmic function takes real values for positive variables. It is worth noting that the boundary values of $\phi_0(z)$ on the upper and lower edges of $(-1, 1)$ are purely imaginary.
   Frequent use will also be made of the following  scaled variables and the     meromorphic  functions, namely,
\begin{equation*}
\tau_\alpha=\frac{z-\alpha}{\alpha+1},\quad \tau_\beta=\frac{z-\beta}{\beta-1},\quad f_s(\tau)=\frac{1}{12b\tau}+\frac{5}{48b\tau^2},\quad \mbox{and}~~f_r(\tau)=\frac {5 \sqrt{2/b}}{72\tau^2\xi_0(\tau)}-f_s(\tau),
\end{equation*}where  $ \xi_{0}(\tau)=\sum^\infty_{k=0}   (-1)^k \sqrt{  b /2}\; B(k+2,    1 /2) \tau^k$, with $B(\xi, \eta)$ being the Beta function.

To describe the behavior at the soft edges, we use conformal mappings
\begin{equation*}
 \lambda_\alpha(z)=e^{-4\pi i/3} \left ( \frac 3 2\right )^{2/3} n^{1/3} \phi^{2/3}_\alpha (z),\quad \mbox{and} \quad \lambda_\beta(z)= \left ( \frac 3 2\right )^{2/3} n^{1/3} \left (-\phi_\beta (z)\right )^{2/3},
\end{equation*}
respectively in neighborhoods of $\alpha$ and $\beta$, and constructed in terms of the $\phi$-functions defined
in Section \ref{MRS} and
analyzed in detail in Section \ref{phi-alpha-phi-beta}.  Briefly, the function $\phi_\alpha(z)$ is analytic in $\mathbb{C}\backslash\{(-\infty, -1]\cup[\alpha,\infty)\}$ such that $n\phi_\alpha (z)=ic_\alpha(n)\tau_\alpha^{3/2}(1+O(\tau_\alpha))$ for small $\tau_\alpha$, with $\arg \tau_\alpha\in (0, 2\pi)$
and $c_\alpha(n)\sim \frac {2\sqrt{2b}} 3 \sqrt n$ for large $n$.
While $\phi_\beta(z)$ is analytic in $\mathbb{C}\backslash(-\infty,\beta]$, such that
$n\phi_\beta (z)=-c_\beta(n)\tau_\beta^{3/2}(1+O(\tau_\beta))$, with $\arg \tau_\beta\in (-\pi, \pi)$ for small $\tau_\beta$, again   $c_\beta(n)\sim \frac {2\sqrt{2b}} 3 \sqrt n$ for large $n$.
One easily obtains $\lambda_\alpha(z)\sim -(2b)^{1/3} \tau_\alpha$  and $\lambda_\beta(z)\sim  (2b)^{1/3} \tau_\beta$ for large $n$, respectively in neighborhoods of $\tau_\alpha=0$ and $\tau_\beta=0$.

Now we are in a position to state the uniform asymptotic expansions for the  orthonormal polynomials $p_n(z)$ as $n\rightarrow \infty$.

{\thm{\label{Main-theorem}
For $r\in (0,1)$, the following holds (see Figure \ref{region-asymptotic} for the regions):

\noindent
\begin{enumerate}[(i)]

\item For $z\in A_r$,
\begin{align}\label{thm-pi-ar}
p_n(z)&=\frac{1}{2\sqrt{b\pi}}e^{-n\phi_\beta(z)}D(z)\left (z+\sqrt{z^2-1}\; \right )^{b/\sqrt{z^2-1}}
\sin\left (b\pi/\sqrt{z^2-1}\right )\nonumber\\
&\quad \times \left[\varrho+\varrho^{-1}+\varrho f_s(\tau_\beta)-\varrho^{-1}f_s(\tau_\alpha)\right]
\left(1+O\left(n^{-1/2}\right)\right),
\end{align}
where $D(z)=\exp\left \{\frac{\sqrt{2}}{2}(1-\frac{1}{6b})\sqrt{z-\beta}
(\sqrt{z+1}-\sqrt{z-\alpha})+O(1/\sqrt{n})\right \}$; cf. \eqref{N-expression}, and $\phi_\beta(z)$ is a function analytic in $\mathbb{C}\setminus (-\infty, \beta]$, as  given in \eqref{phi-beta}.
\item For $z\in B_r$,
\begin{align}\label{thm-pi-br}
p_n(z)&=\frac{ (-1)^{n}}{2\sqrt{b\pi}}(z+\sqrt{z^2-1})^{b/\sqrt{z^2-1}}
(1-e^{-2i\pi b/\sqrt{z^2-1}})^{1/2}\nonumber\\
&\quad\times\left[\left(1+f_s(\tau_\beta)+O\left(n^{-1/2}\right)\right)\varrho e^{-\pi i/4}\cos \Theta_B\right.\nonumber\\
&\quad\left.+\left(1-f_s(\tau_\alpha)+O\left(n^{-1/2}\right)\right)\varrho^{-1} e^{\pi i/4}\sin \Theta_B\right],
\end{align}
where $\Theta_B=i n\phi_\alpha(z)+(i\ln \phi_0(z))/2-i\ln D(z)+\pi/4$ with $D(z)$ the same as that in \eqref{thm-pi-ar}, $\phi_0(z)$ is defined in \eqref{phi-0} and $\phi_\alpha(z)$ is defined as \eqref{phi-alpha}.

\item For $z\in C_r$,
\begin{align}\label{thm-pi-cr}
p_n(z)&=\frac{e^{b}}{2\sqrt{b\pi}}
\left[
\left(1+f_s(\tau_\beta)+O\left(n^{-1/2}\right)\right) \varrho e^{-\pi i/4}\cos \Theta_C
\right.\nonumber\\
&~~~\left.
+\left(1+O\left(n^{-1/2}\right)\right)\varrho^{-1} e^{\pi i/4}\sin \Theta_C\right],
\end{align}
where $\Theta_C=i n\phi_\beta(z)+b\pi/\sqrt{z^2-1}-\pi/4+O(1/n)$.
\begin{figure}[t]
\begin{center}
\includegraphics[width=12cm,bb=37 362 546 510]{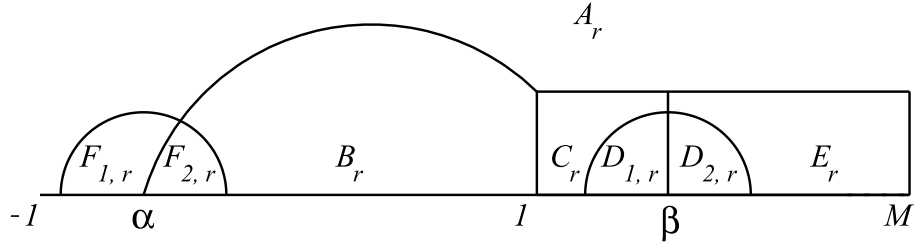}
 \caption{Regions of uniform asymptotic approximations, where $M$ is a real constant such that $M>x_1=\sqrt{1+b^2}$.}
 \label{region-asymptotic}
\end{center}
\end{figure}

\item For $z\in D_{1,r}\cup  D_{2, r}~(|\tau_\beta|<r)$,
    \begin{align}\label{pi-n-beta}
       p_n(z)&=\frac{e^{b}}{2\sqrt{b}}\left\{n^{1/12} A_1(z)\{\lambda_\beta(z)\}^{1/4}\varrho^{-1}\left[ 1+O\left(n^{-1/2}\right)\right]\right.\nonumber\\
        &~~~\left.+n^{-1/12} A_2(z)\{\lambda_\beta(z)\}^{- 1/4}\varrho
        \left[1-f_r(\tau_\beta)+O\left(n^{-1/2}\right)\right]\right\},
    \end{align}
         where  the pair of functions
    $$
     A_1(z)=-\cos\left (b\pi/\sqrt{z^2-1}\right )\Ai(n^{1/3}\lambda_\beta(z))+\sin\left (b\pi/\sqrt{z^2-1}\right )\Bi(n^{1/3}\lambda_\beta(z)),$$
     $$    A_2(z)=-\cos\left (b\pi/\sqrt{z^2-1}\right )\Ai'(n^{1/3}\lambda_\beta(z))+\sin\left (b\pi/\sqrt{z^2-1}\right )\Bi'(n^{1/3}\lambda_\beta(z)).
    $$
\item For $z\in E_r$,
\begin{align}\label{thm-pi-er}
&p_n(z)= \frac{1}{2\sqrt{b\pi}}{(z+\sqrt{z^2-1})}^{b/\sqrt{z^2-1}}
\left\{\left[\varrho+\varrho^{-1}+\varrho
f_s(\tau_\beta)-\varrho^{-1}f_s(\tau_\alpha)\right]\right.\nonumber\\
&\left.\times\sin(b\pi/\sqrt{z^2-1})e^{-n\phi_\beta(z)}
\left(1+O\left(n^{-1/2}\right)\right)+O(n^{1/4}e^{\Re(n\phi_\beta(z)-ib\pi/\sqrt{z^2-1})})\right\}.
\end{align}

\item For $z\in F_{1, r}\cup  F_{2, r}~(|\tau_\alpha|<r)$,
    \begin{align}\label{pi-n-alpha}
        p_n(z)&=\frac{(-1)^n}{2\sqrt{b}}e^{-b+\frac{b\pi}{\sqrt{2(z+1)}}}
        \left\{n^{1/12} \{\lambda_\alpha(z)\}^{1/4}\varrho\; \Ai(n^{1/3}\lambda_\alpha(z))\left[1+O\left(n^{-1/2}\right)\right] \right.\nonumber\\
        & \quad\left.-        n^{-1/12} \{\lambda_\alpha(z)\}^{-1/4}
        \varrho^{-1}\Ai'(n^{1/3}\lambda_\alpha(z))\left[2-\frac{1}{6b} + f_r(\tau_\alpha)
        +O\left(n^{-1/2}\right)\right]\right\}.
    \end{align}

\end{enumerate}}}

Later in Section \ref{comparison}, the results will be compared with those of Wang and Zhao \cite{wang-zhao},
obtained earlier via integral methods.

\section{ RHP and IP for the sieved Pollaczek polynomials}
\setcounter{equation} {0}  \label{RHP-IP}
We begin with the following mixed RHP and IP formulation. The formulation has been given in Wang \cite{wang}. For an earlier RHP version, see \cite{Fokas}, while an IP version can be found in \cite{Borodin}. The  composite RHP and IP is as follows:

\begin{description}
  \item[($Y_a$)]     $Y(z)$ is analytic in
  $\mathbb{C}\backslash([-1,1]\cup  \mathcal{X} )$, $\mathcal{X}=\left \{x_k: \;  x_k=\sqrt{1+\frac{b^2}{k^2}},~k=1,2,\cdots\right\}$.

  \item[($Y_b$)]    $Y(z)$  satisfies the jump condition
\begin{equation}\label{Y-jump-continuous}
  Y_+(x)=Y_-(x) \left(
                               \begin{array}{cc}
                                 1 & w_c(x) \\
                                 0 & 1 \\
                                 \end{array}
                             \right),
 ~~~ ~~x\in (-1,1),\end{equation}where $w_c(x)$ is  the absolutely continuous weight defined in \eqref{weight-continuos}
for   $x\in (-1, 1)$.
\item[($Y_c$)]    $Y(z)$ has simple poles at the nodes $x_k=\sqrt{1+\frac{b^2}{k^2}}$,  and  satisfies the residue  condition
\begin{equation}\label{Y-residue}
  {\Res}_{z=x_k}Y(z)=\lim_{z\rightarrow x_k}Y(z) \left(
                               \begin{array}{cc}
                                 0& -w_d(x_k)/2\pi i \\
                                 0 & 0 \\
                                 \end{array}
                             \right)
~~\mbox{for}~~k=1,2,\cdots,\end{equation}
where the discrete  weight is
\begin{equation}\label{weight-discontinuous}
   w_d(x)=4x^{-1}(x^2-1)^{3/2}(x-\sqrt{x^2-1})^{\frac{2b}{\sqrt{x^2-1}}},
\end{equation}being positive for $x\in (1, +\infty)$; cf. \eqref{discontinuous-weight}.
We note that $w_d(x_k)\sim 4b^3 e^{-2b}/k^3$ as $k\to\infty$.

\item[($Y_d$)]    The asymptotic behavior of $Y(z)$  at infinity
  is
  \begin{equation}\label{Y-infty}Y(z)=\left (I+O\left (\frac 1 z\right )\right )\left(
                               \begin{array}{cc}
                                 z^n & 0 \\
                                 0 & z^{-n} \\
                               \end{array}
                             \right),~~~\mbox{as}~~z\rightarrow
                             \infty  .\end{equation}
\item[($Y_e$)] $Y(z)$ has the following behavior at $\pm1$.
\begin{equation}\label{Y-pm1}Y(z)= \left(
                               \begin{array}{cc}
                                 O(1) & O(\ln |z\mp 1|) \\
                                 O(1) & O(\ln |z\mp 1|) \\
                               \end{array}
                             \right),~~~\mbox{as}~~z\to \pm 1
                             .\end{equation}

\end{description}
It is readily verified that the unique solution to the RHP for Y is
\begin{equation}\label{Y-def}
  Y(z)=\left(
                               \begin{array}{cc}
                                 \pi_n(z) & \displaystyle{ \frac{1}{2\pi i}\int_{-1}^1\frac{\pi_n(x)w_c(x)dx}{x-z} +\frac{1}{2\pi i} \sum_{k=1}^{\infty}\frac{\pi_n(x_k)w_d(x_k)}{x_k-z}} \\[.4cm]
                                 -2\pi i \gamma^2_{n-1} \pi_{n-1}(z) &\displaystyle{ -\gamma^2_{n-1}   \int_{-1}^1\frac{\pi_{n-1}(x)w_c(x)dx}{x-z} -\gamma^2_{n-1} \sum_{k=1}^{\infty}\frac{\pi_{n-1}(x_k)w_d(x_k)}{x_k-z}
                                 }\\
                                 \end{array}
                             \right)
,\end{equation}
where $\pi_n(z)$ are the monic sieved Pollazeck polynomials, and $p_n(z)=\gamma_n \pi_n(z)$ are the corresponding orthonormal polynomials.

\section{Removing the poles of $Y(z)$: RHP for $U(z)$}
\setcounter{equation} {0}  \label{removing-poles}

We use the notion of band and saturated regions of  Baik {\it et al.}\;\cite{bkmm} to treat the present composite weight. In later sections, we will show that the    real interval $(-1, \alpha)$ is the   void, $(\beta, x_1)$ is the  saturated region, and  $(\alpha, \beta)$ is the   band: part of it belongs to the absolutely continuous support, the other part corresponds to the accumulating nodes. It will be shown that $\alpha\approx -1+b/n$ and $\beta\approx 1+b/n$ for large polynomial degree $n$.

Applying the ideas of Bleher and Liechty \cite{BleherE, Bleher}, we may define
\begin{equation}\label{D-upper}
 D^u_{\pm}(z)=
\left(
  \begin{array}{cc}
  1 &   \frac {\gamma'(z)  w_d(z) e^{\mp  i\pi \gamma(z)}} {2  i\sin (\pi \gamma(z))}  \\
  0 & 1 \\
  \end{array}
\right),~~\pm \Im z \geq 0,
\end{equation}
\begin{equation}\label{D-lower}
 D^l_{\pm}(z)=\pm
\left(
  \begin{array}{cc}
\frac {e^{\pm i\pi \gamma(z)}}{2i \sin(\pi\gamma(z))}  &  0  \\
  \frac {1} {\gamma'(z)w_d(z)} & \frac {2i \sin(\pi\gamma(z))}{e^{\pm i\pi \gamma(z)}} \\
  \end{array}
\right),~~\pm \Im z\geq 0,
\end{equation}where $\gamma(z)=\frac {b}{\sqrt{z^2-1}}$ and $w_d(z)= 4z^{-1}(z^2-1)^{3/2}(z-\sqrt{z^2-1})^{\frac{2b}{\sqrt{z^2-1}}}$, branches are chosen such that  $\gamma(z)$ and $w_d(z)$ are analytic respectively  in $\mathbb{C}\backslash [-1, 1]$ and  $\mathbb{C}\backslash (-\infty, 1]$, and both functions are  positive for $z\in (1, +\infty)$.

As in Bleher and Liechty \cite{Bleher}, we introduce
\begin{equation}\label{U-upper}
U^u(z)=Y(z)\left\{ \begin{array}{ll}
                    D^u_+(z),  & \Im z\geq 0, \\
                    D^u_-(z),  & \Im z\leq 0,
                   \end{array}\right .
\end{equation}
and
\begin{equation}\label{U-lower}
U^l(z)=Y(z)\left\{ \begin{array}{ll}
                    D^l_+(z),  & \Im z\geq 0, \\
                    D^l_-(z),  & \Im z\leq 0.
                   \end{array}\right .
\end{equation}
It is readily verified that   both functions $U^l(z)$ and $U^u(z)$ are analytic in $\{z|\Re z>1,~\Im z\geq 0\}$, and in $\{z|\Re z>1,~\Im z\leq 0\}$, that is, all simple poles of $Y(z)$ at $z=x_k=\sqrt{1+b^2/k^2}$, $k=1,2,\cdots$, have been removed.

We introduce several scalar  auxiliary functions.
\begin{equation}\label{d-E}
 d_E(z)=\frac {\Gamma(1-\gamma(z))} {\sqrt {2\pi}} e^{ -\gamma(z)+\left (\gamma(z)-\frac 1 2\right )\ln (-\gamma(z))},~~z\in \mathbb{C}\setminus [-1, \infty),
\end{equation}
\begin{equation}\label{d-I}
  d_I(z)=\frac {\sqrt{2\pi}} { \Gamma(\gamma(z))}   e^{ -\gamma(z)+\left (\gamma(z)-\frac 1 2\right )\ln  \gamma(z)},~~z\in \mathbb{C}\setminus (-\infty, 1],
\end{equation}and
\begin{equation}\label{chi}
  \chi (z)=2\sqrt{\pi b} \left (e^{\mp \pi i} \varphi(z)\right )^{-\frac 1 2}~~\mbox{as}~~\pm \Im z>0,~~z\in \mathbb{C}\setminus [-1, \infty),
\end{equation}where $\varphi(z)=z+\sqrt{z^2-1}$ is analytic in $\mathbb{C}\setminus [-1, 1]$ such that $\varphi(z)  \approx 2z$ for large $z$,
   and the branches of the  logarithms
are chosen such that $$\arg(-\gamma(z)) =\pm \pi+\arg(\gamma(z)) =\pm \pi -\frac 1 2 \left (\arg(z-1)+\arg(z+1)\right ), ~~\pm \Im z>0.$$

We note that $d_E(z)$ and  $d_I(z)$ approaches $1$ as $z\to 1$ respectively from $\Re z\leq 1$ and $\Re z \geq 1$, and that
 $d_E(z)\chi(z)=1+O( \ln |z|/z)$ as $z\to\infty$.  The functions  $d_E(z)$ and
 $d_I(z)$ will simplify the jumps,    and $\chi (z)$ will normalize  the behavior at infinity,  of the RHPs for $U$, $T$ and $S$ in later sections,
while   the behavior   at $z=1$ retains the original form. Similar auxiliary functions have been used
 in Wang and Wong \cite{WangWong}  and Lin and Wong \cite{LinWong}.

Define a matrix-valued function $U(z)$ as
\begin{equation}\label{U-def}
U(z)=\left\{
\begin{array}{ll}
  U^l(z)\left \{ d_I(z) \chi(z)\right \}^{\sigma_3},  &  z\in (\beta, M)\times (0, \pm i\varepsilon), \\[.2cm]
  U^u(z)\left \{ d_I(z) \chi(z)\right\}^{\sigma_3} ,  &  z\in (1, \beta)\times (0, \pm i\varepsilon),\\[.2cm]
  Y(z)\left \{ d_E(z) \chi(z)\right \}^{\sigma_3}, &\mbox{otherwise},
\end{array}
\right .
\end{equation}
where $M$ is a constant such that $M>x_1=\sqrt{1+b^2}$, $\varepsilon=\delta/\sqrt n$ for a small positive $\delta$ independent of  $z$, and $\beta=\beta_n$  in \eqref{U-def} is one of the MRS numbers to be determined. It will be shown later that    $\beta\approx 1+\frac b n$ for large $n$.

\begin{figure}[t]
\begin{center}
\includegraphics[width=12cm,bb=31 346 570 442]{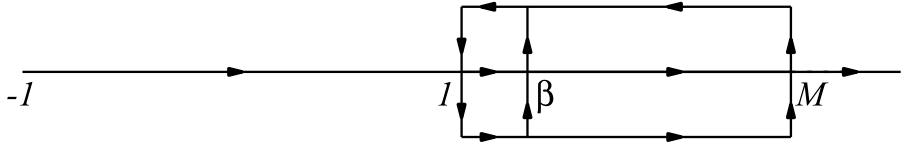}
 \caption{The contour $\Sigma_U$ for $U(z)$.}
 \label{contour-for-U}
\end{center}
\end{figure}

It is readily verified that $U(z)$ solves the following RHP:

\begin{description}
  \item[($U_a$)]     $U(z)$ is analytic in
  $\mathbb{C}\backslash \Sigma_U$, where $\Sigma_U$ is illustrated in Figure \ref{contour-for-U}.

  \item[($U_b$)]    $U(z)$  satisfies the jump condition
\begin{equation}\label{U-jump}
  U_+(z)=U_-(z) J_U(z),
 ~~~ z\in \Sigma_U,\end{equation}where the jump on the real axis is
\begin{equation}\label{U-jump-real-axis}
J_U(x)= \left\{
\begin{array}{ll}
 \left(
          \begin{array}{cc}
           \left (\phi_0\right )_+(x)   & r(x) \\
            0 & \left (\phi_0\right )_-(x) \\
          \end{array}
        \right)
  ,  &  x\in (-1, 1 ), \\[.5cm]
 \left(
          \begin{array}{cc}
           -1  & r(x) \\
            0 & -1 \\
          \end{array}
        \right)
  ,  &  x\in (1, \beta ), \\[.5cm]
     \left(
        \begin{array}{cc}
          e^{2\pi i\gamma(x)}  & 0 \\
           \frac {4\pi b d_I^2(x) } {\varphi(x)   \gamma'(x) w_d(x)} &  e^{-2\pi i\gamma(x)} \\
        \end{array}
      \right)
  , & x\in (\beta, M),\\[.5cm]
\left(
        \begin{array}{cc}
          e^{2\pi i\gamma(x)}  & 0 \\
           0 &  e^{-2\pi i\gamma(x)} \\
        \end{array}
      \right)
  , & x\in (M, \infty ),
\end{array}
\right .
\end{equation}
where $\phi_0(z)$ is the scalar function defined  in \eqref{phi-0}.
 We note that for $x\in (-1, 1 )$, both $\left |\left (\phi_0\right )_\pm(x)\right |=1$, and such that $\left (\phi_0\right )_+(x) \left (\phi_0\right )_-(x)=1$,
  and
\begin{equation}\label{r-def}
r(x)= \left\{
 \begin{array}{ll}
  \frac {w_c(x)} {4\pi b} \left (1- \exp\left ( -\frac {2\pi b}{\sqrt{1-x^2}}\right )\right )=\frac 1 \pi \exp\left ( -\frac {2b\arccos x} {\sqrt{1-x^2}}\right ),  & x\in (-1, 1), \\[.4cm]
  \frac {-\gamma'(x) w_d(x)} {4\pi b} \frac {\varphi(x)} {  d_I^2(x) },  &x\in (1, \beta)
 \end{array}
\right .
\end{equation}
can  actually be extended to  a  continuous function in $(-1, \infty)$. While the jumps on the off-real-axis contours are
\begin{equation}\label{U-jump-off-real-axis}
J_U(z)= \left\{
\begin{array}{ll}
 \left(
          \begin{array}{cc}
           1  & 0 \\
          \pm\frac {d_E(z) d_I(z) \chi^2(z)} { \gamma'(z) w_d(z)}  & 1 \\
          \end{array}
        \right)
  ,  &  \begin{array}{l}
                                                               z\in (\beta, M)\pm i\varepsilon,\\
                                                              z\in M\pm i(0,\varepsilon) ,
                                                             \end{array}
  \\[.5cm]
\left(
          \begin{array}{cc}
          \frac {\pm 2i\sin (\pi\gamma(z))}  {e^{\pm \pi i \gamma(z)} }   & \frac  {\gamma'(z) w_d(z) e^{\mp \pi i \gamma(z)} }   {2i\sin (\pi\gamma(z))   d_E(z) d_I(z) \chi^2(z)}  \\
         0 & \frac {e^{\pm \pi i \gamma(z)} } {\pm 2i\sin (\pi\gamma(z))}\\
          \end{array}
        \right)
  ,
   &\begin{array}{l}
              z\in (1, \beta)\pm i\varepsilon ,     \\
              z\in 1\pm i(0,\varepsilon) ,
               \end{array}
\\[.5cm]
\left(
  \begin{array}{cc}
    \frac {\pm 2i \sin\pi\gamma(z)} {e^{\pm i\pi\gamma(z)}} &  \frac {\pm \gamma'(z) w_d(z)e^{\mp 2 i\pi\gamma(z)} }{d_I^2(z) \chi^2(z)}  \\
    \frac {d_I^2(z) \chi^2(z)} {\mp  \gamma'(z) w_d(z)  } & 1 \\
  \end{array}
\right)
,
   & \begin{array}{l}
   z\in \beta \pm i(0,\varepsilon).
     \end{array}

\end{array}
\right .
\end{equation}

\item[($U_c$)]      The asymptotic behavior of $U(z)$  at infinity
  is
  \begin{equation}\label{U-infty}U(z)=\left (I+O\left (\frac {\ln z} z\right )\right )\left(
                               \begin{array}{cc}
                                 z^n & 0 \\
                                 0 & z^{-n} \\
                               \end{array}
                             \right),~~~\mbox{as}~~z\rightarrow
                             \infty  .\end{equation}
\item[($U_d$)]
   $U(z)$ has the following behavior
  \begin{equation}\label{U-pm1}
 U(z)=O\left(\ln|z\mp 1|\right),\quad \mbox{as}~z\rightarrow \pm 1.
  \end{equation}

\end{description}

\section{MRS numbers and auxiliary functions}
 \setcounter{equation} {0}  \label{MRS}
Assume that $\psi(x) dx$ is the equilibrium measure, supported on $(\alpha, \infty)$. We consider the following $g$-function, to be used in the transformation \eqref{transformation-U-to-T} below.
\begin{equation}\label{g-def}
g(z)= \int_\alpha^\infty \ln (z-x) \psi(x) dx, ~~~z\in \mathbb{C}\backslash \mathbb{R},
\end{equation}in which the branch of the logarithm is chosen such that $\arg (z-x)\in (-\pi, \pi)$,  $\psi(x)$ is to be determined for $x\in (\alpha, \beta)$, and we take
\begin{equation}\label{psi-(beta-infty)-def}
\psi(x)=-\frac {\gamma'(x)} n =\frac 1 n  bx(x^2-1)^{-3/2}~~\mbox{for}~~  x\in  (\beta, \infty),
\end{equation}
understanding that
 $-\gamma'(x)=bx(x^2-1)^{-3/2}$ is the limit density of the nodes $x_k=\sqrt{1+\frac {b^2}{k^2}}$, $k=1,2,\cdots$.

The $g$-function can be determined by a phase condition of the form
\begin{equation}\label{g-phase-condition}
g_+(x)+g_-(x)-l +\frac 1 n  \ln r(x)  =0, ~~~x\in (\alpha, \beta),
\end{equation} where $l$ is the Lagrange multiplier independent of $x$.

It is readily seen that $r'(x)/r(x)$ is an infinitely smooth function in $(-1, 1)\cup(1, \infty)$, such that
\begin{equation}\label{r-end-behavior}
\frac {r'(x)}{r(x)}=\left \{
\begin{array}{ll}
   \frac {b\pi}{\sqrt 2 (1+x)^{3/2}} +  O\left ( \frac 1 {\sqrt{1+x}}  \right ) & \mbox{as}~~x\to -1^+, \\[.2cm]
 \frac {2b} 3 + O(x-1) & \mbox{as}~~x\to 1^-, \\[.2cm]
 O\left ( \frac 1 {\sqrt{x-1}}\right ) & \mbox{as}~~x\to 1^+ .
\end{array}\right .
\end{equation}
Denoting
\begin{equation}\label{G-def}
G(z)= \frac 1 {\pi i} \int^\infty_\alpha \frac { \psi(x) dx} {x-z}~~\mbox{for}~~z\in \mathbb{C}\backslash [\alpha, \infty),
\end{equation}  we see that  $G(z)=-\frac 1 {\pi i} g'(z)$ in their common domain of analyticity, and that $G(z)$ solves the scalar RHP
\begin{equation}\label{G-jump}
\left\{
\begin{aligned}
    G_+(x)+G_-(x)&= \frac 1 {n\pi i}\frac {r'(x)}{r(x)}~&\mbox{for}~x\in (\alpha, \beta),\\
    G_+(x)-G_-(x)&=\frac 2 n bx(x^2-1)^{-3/2}~&\mbox{for}~x\in (\beta, \infty)
\end{aligned}
\right .\end{equation}
 and
\begin{equation}\label{G-infty} G(z)=-\frac 1 {\pi i} \frac 1 z+O\left (\frac {\ln z} {z^2}\right )~~\mbox{as}~z\to\infty.
\end{equation}
Such a solution   can be expressed as
\begin{equation}\label{G-solution}
G(z) =  \frac {\sqrt{(z-\alpha)(z-\beta)}} {2 n\pi i} \left [  \int^\beta_\alpha \frac  {-r'(x)/r(x)}{\pi\sqrt{(x-\alpha)(\beta-x)}} \frac {dx}{x-z}+\int^\infty_\beta \frac {2bx(x^2-1)^{-3/2}}{ \sqrt{(x-\alpha)(x-\beta)}} \frac {dx}{x-z}\right ],
\end{equation}subject to the condition \eqref{G-infty} at infinity, which now takes the form
\begin{equation}\label{MRS-numbers-eqn}
\left\{\begin{aligned}
\int^\beta_\alpha \frac  {r'(x)/r(x)\; dx}{\pi\sqrt{(x-\alpha)(\beta-x)}} -\int^\infty_\beta \frac {2bx(x^2-1)^{-3/2} dx }{ \sqrt{(x-\alpha)(x-\beta )}} & =0,\\
\int^\beta_\alpha \frac  {x r'(x)/r(x)\; dx}{\pi\sqrt{(x-\alpha)(\beta-x)}} -\int^\infty_\beta \frac {2bx^2(x^2-1)^{-3/2} dx }{ \sqrt{(x-\alpha)(x-\beta)}} &=-2  n.
\end{aligned}\right .
\end{equation}

To determine the MRS numbers $\alpha$ and $\beta$, we include some details in what follows.

In view of \eqref{r-end-behavior}, from \eqref{MRS-numbers-eqn} we have
$$\left\{\begin{aligned}
\frac {b } 2 \int^\beta_\alpha \frac  {dx}{ (x+1)^{3/2} \sqrt{x-\alpha}} -\frac {b } 2 \int^\infty_\beta \frac {dx }{(x-1)^{3/2} \sqrt{x-\beta}} & =O(1),\\
-\frac {b } 2 \int^\beta_\alpha \frac  {dx}{ (x+1)^{3/2} \sqrt{x-\alpha}} -\frac {b } 2 \int^\infty_\beta \frac {dx }{(x-1)^{3/2} \sqrt{x-\beta}} & = -2  n+O(1).
\end{aligned}\right .$$
Noticing that the indefinite integrals
$$\int \frac  {dx}{ (x+1)^{3/2} \sqrt{x-\alpha}}=\frac 2 {\alpha+ 1} \frac {\sqrt{x-\alpha}} {\sqrt{x+1}}~~\mbox{and}~~\int  \frac {dx }{(x-1)^{3/2} \sqrt{x-\beta}}=\frac 2 {\beta-1} \frac {\sqrt{x-\beta}} {\sqrt{x-1}} $$up to  an arbitrary constant,
  the equations in \eqref{MRS-numbers-eqn} can be written in the form
$$
\frac {b}{\alpha+ 1}-\frac {b}{\beta-1} =O(1),~~\mbox{and}~~  -\frac {b}{\alpha+1}-\frac { b}{\beta-1} = -2  n+O(1)
$$ for large $n$. Hence we have
\begin{equation}\label{MRS-numbers-approx}
 \alpha =-1+\frac b n+O\left(\frac{1}{n^2}\right),~~\mbox{and}~~  \beta  = 1+\frac b n+O\left(\frac{1}{n^2}\right),~~\mbox{as}~n\to\infty. \end{equation}Refinements can be obtained by using a similar argument as in \cite{zhou-zhao}, or in \cite{zhou-xu-zhao}.

The equilibrium measure can be expressed as
\begin{equation}\label{equilibrium-measure}
\psi(x)=\Re G_+(x)~~\mbox{for}~~x\in (\alpha, \beta).\end{equation}

We proceed to define certain $\phi$-functions. To this aim, we seek a function $\nu_\alpha(z)$, analytic in $\mathbb{C}\backslash \{ (-\infty, -1]\cup [\alpha, +\infty)\}$, such that
\begin{equation}\label{nu-boundary}
\left (\nu_\alpha\right )_\pm(x)=\pm \pi i \psi(x),~~x\in (\alpha, 1).\end{equation}
From \eqref{G-def}, and using the Sokhotski-Plemelj formula, we have
$$G_\pm(x)=\pm \psi(x) +\frac 1 {\pi i} p.v. \int_\alpha^{+\infty}\frac {\psi(t) dt} {t-x},~~x\in (\alpha, \infty),$$
 in which the boundary values share the same imaginary part. Further more, in view of \eqref{G-jump}, it is readily seen that
 $$\Im G_+(x)=\Im G_-(x)=-\frac 1 {2n\pi} \frac {r'(x)}{r(x)} = \frac b { n\pi} \frac {x\arccos x-\sqrt{1-x^2}} {(1-x^2)^{3/2}}~~\mbox{for}~~x\in (\alpha, 1);$$ cf. \eqref{weight-continuos}  and \eqref{r-def}.
 It is worth noting that $G(z)$ is analytic in $\mathbb{C}\setminus [\alpha, \infty)$, and   $-\frac {r'(x)}{r(x)}$ can be analytically extended from  $x\in (-1, 1)$ to
\begin{equation}\label{h-alpha}
h_\alpha(z)=  -  \frac {2b\left [\sqrt{1-z^2} +iz \ln   (z+i\sqrt{1-z^2}\;  )\right ] } {(1-z^2)^{ 3/ 2}}~~\mbox{for}~z\in \mathbb{C}\setminus (-\infty, -1],
\end{equation}
where  branches are chosen such that $\arg (z+1)\in (-\pi, \pi)$ and $\arg(1-z)\in (-\pi, \pi)$.     Therefore, we chose
\begin{equation}\label{nu-representation}
\nu_\alpha(z)=\pi i G(z) -\frac b {n} \frac {\sqrt{1-z^2} +iz \ln   (z+i\sqrt{1-z^2}\;  )} {(1-z^2)^{ 3/ 2}},~~z\in \mathbb{C}\backslash \{ (-\infty, -1]\cup [\alpha, +\infty)\}.\end{equation}
Consequently, we can define a $\phi$-function
\begin{equation}\label{phi-alpha}
\phi_\alpha(z)=\int^z_\alpha \nu_\alpha(\zeta) d\zeta~~\mbox{for}~~z\in \mathbb{C}\backslash \{ (-\infty, -1]\cup [\alpha, +\infty)\}.\end{equation}

Now we turn to the other critical point $x=\beta\approx 1+\frac b n$. Similar to the above derivation, we  introduce
\begin{equation}\label{nu-beta}
\nu_\beta(z)=\pi i \left \{ G(z) \pm \frac 1 n \gamma'(z)\right\}  -\frac 1 {2n} \frac {r'(z)}{r(z)}~~\mbox{for}~~\pm \Im z >0, \end{equation}where
$r(z)= \frac {-\gamma'(z) w_d(z)} {4\pi b} \frac {\varphi(z)} {  d_I^2(z) }$; cf. \eqref{r-def}, with
 $ -\gamma'(z)w_d(z)=4b\left ( z-\sqrt{z^2-1}\right )^{\frac {2b} {\sqrt{z^2-1}}} $   being  analytically extended to $\mathbb{C}\setminus (-\infty, -1]$,   branches chosen such that
$\arg (z\pm 1)\in (-\pi, \pi)$, and real positive  for $z=x\in (1, +\infty)$.

Therefore, $r(z)$, and hence $r'(z)/r(z)$, are analytic in $\mathbb{C}\setminus (-\infty, 1]$. While $\gamma(z)=\frac b {\sqrt{z^2-1}}$ is analytic in $\mathbb{C}\setminus [-1, 1]$, and positive for $z>1$.
It is readily seen that $\nu_\beta(z)$ is analytic in  $\mathbb{C}\setminus (-\infty, \beta]$, such that
$$ \left (\nu_\beta\right )_\pm (x) = \pm \pi i\left \{ \psi(x)+\frac {1} n \gamma'(x)\right\}~~\mbox{for}~~x\in (1, \beta).$$ Here, use has been made of the phase condition \eqref{g-phase-condition}. Upon these we define another $\phi$-function
\begin{equation}\label{phi-beta}
\phi_\beta(z)=\int^z_\beta \nu_\beta (\zeta) d\zeta~~\mbox{for}~~z\in \mathbb{C}\backslash   (-\infty,   \beta] .\end{equation}

We derive connections between these auxiliary functions. Indeed,  substituting  \eqref{nu-representation} into \eqref{phi-alpha} yields
\begin{equation}\label{phi-alpha-g-relation}
 \phi_\alpha (z)+g(z)+\frac {i b \ln(z+i\sqrt{1-z^2}\; )}  {n\sqrt{1-z^2}}  -\frac l 2-\frac {\ln\pi}{2n} \mp \pi i=0,~~\pm \Im z >0.
  \end{equation}Here use has been made of the representation \eqref{g-def} and the phase condition \eqref{g-phase-condition}. The constant $l$ is the same as in the phase condition \eqref{g-phase-condition}, and can be determined  from \eqref{phi-alpha-g-relation} as
  $$l=2\int^\infty_\alpha \ln (x-\alpha)\psi(x) dx- \frac {2b\arccos \alpha }{ n\sqrt{1-\alpha^2}}-\frac {\ln\pi}{n};$$ see \eqref{l-asymptotic} for an asymptotic approximation of $l$.

Similarly, we have
\begin{equation}\label{phi-beta-g-relation}
 \phi_\beta (z)+g(z)\mp \frac {\pi i}  n  \gamma(z) +\frac 1 {2n} \ln r(z)   -\frac l 2=0,~~\pm \Im z >0
  \end{equation}for the same $l$, where $r(z)$ is as defined in  \eqref{nu-beta},  and $ \ln r(z) $ is analytic in $\mathbb{C}\backslash (-\infty, 1]$ such that $\arg r(z) =0$ for $z=x\in (1, +\infty)$.

\section{Nonlinear steepest descent analysis}
\setcounter{equation} {0}  \label{steepest-descent-analysis}
The procedure consists of a series of invertible transformations $U(z)\to T(z)\to S(z)\to R(z)$.
The pioneering work in this respect includes  Deift and Zhou {\it{et al.}}
\cite{DeiftZhou,dkmv2}; see also  Bleher and  Its \cite{bleher-its}.
 To accomplish the transformations, auxiliary functions are analyzed and parametrices are constructed at critical points $z=\alpha$ and $z=\beta$, as well as in a region around the infinity.

\subsection{The first transformation $U\to T$}
The first transformation $U\to T$ is the following  normalization at infinity
\begin{equation}\label{transformation-U-to-T}
U(z)= e^{\frac 1 2 nl\sigma_3} T(z) e^{n\left (g(z)-\frac 1 2 l\right ) \sigma_3}, ~~z\in \mathbb{C}\backslash \Sigma_U,
\end{equation}
where $l$ is a constant; cf. \eqref{g-phase-condition}, $\Sigma_U$ is illustrated in Figure \ref{contour-for-U}, $\sigma_3$ is one of  the Pauli matrices,  defined as
\begin{equation}\label{Pauli-matrices}
\sigma_1=\left(
                   \begin{array}{cc}
                     0 &1 \\
                    1 & 0 \\
                   \end{array}
                 \right),   ~~\sigma_2=\left(
                   \begin{array}{cc}
                     0 & -i \\
                    i & 0 \\
                   \end{array}
                 \right),~~\mbox{and}~\sigma_3=\left(
                   \begin{array}{cc}
                     1 & 0 \\
                    0 & -1 \\
                   \end{array}
                 \right).
 \end{equation}

It is readily verified that   $T(z)$
solves the following RHP
\begin{description}
  \item[($T_a$)]     $T(z)$ is analytic in
  $\mathbb{C}\backslash \Sigma_U$, with  $\Sigma_U$   illustrated in Figure \ref{contour-for-U}.

  \item[($T_b$)]    $T(z)$  satisfies the jump condition
\begin{equation}\label{T-jump}
  T_+(z)=T_-(z) J_T(z),
 ~~~ z\in \Sigma_U,\end{equation}where the jumps are
 \begin{equation}\label{T-jump-on-R}
    J_T(x)= \left\{
\begin{array}{ll}
 \left(
          \begin{array}{cc}
           \frac 1 { \left(\phi_0\right ) _-(x)} & e^{-2n\phi_\alpha(x)}  \\
            0 & \frac 1 {\left(\phi_0\right ) _+(x)}\\
          \end{array}
        \right)
  ,  &  x\in (-1, \alpha ), \\[.6cm]
     \left(
        \begin{array}{cc}
         \frac { e^{-2n\left(\phi_\alpha\right ) _-(x)}}{ \left(\phi_0\right ) _-(x)}  & 1 \\
           0 & \frac {e^{-2n\left(\phi_\alpha\right ) _+(x)}}{ \left(\phi_0\right ) _+(x)}\\
        \end{array}
      \right)
  , & x\in (\alpha, 1),\\[.6cm]
  \left(
        \begin{array}{cc}
          -e^{-2\pi i \gamma(x)-2n\left(\phi_\beta\right ) _-(x) }  & 1 \\
           0 & -e^{2\pi i \gamma(x)-2n\left(\phi_\beta\right ) _+(x)  }\\
        \end{array}
      \right)
  , & x\in (1, \beta),\\[.6cm]
   \left(
        \begin{array}{cc}
         1  & 0 \\
         - e^{2n\phi_\beta(x)} &  1\\
        \end{array}
      \right)
  , & x\in (\beta, M),
    \end{array}
 \right .
\end{equation}
 and
 \begin{equation}\label{T-jump-off-real-axis}
J_T(z)= \left\{
\begin{array}{ll}
 \left(
          \begin{array}{cc}
           1  & 0 \\
           \frac {e^{2n\phi_\beta (z)\mp \pi i \gamma(z)}} {2 i\sin \pi\gamma(z)}  & 1 \\
          \end{array}
        \right)
  ,  &   \begin{array}{l}
                                                               z\in (\beta, M)\pm i\varepsilon, \\
                                                              z\in M\pm i(0,\varepsilon) ,
                                                             \end{array}
  \\[.5cm]
\left(
          \begin{array}{cc}
          \frac {\pm 2i\sin (\pi\gamma(z))}  {e^{\pm \pi i \gamma(z)} }   & \pm e^{-2n\phi_\beta (z)}  \\
         0 & \frac {e^{\pm \pi i \gamma(z)} } {\pm 2i\sin (\pi\gamma(z))}\\
          \end{array}
        \right)
  ,
   &\begin{array}{l}
              z\in (1, \beta)\pm i\varepsilon ,     \\
              z\in 1\pm i(0,\varepsilon) ,
               \end{array}
\\[.5cm]
\left(
  \begin{array}{cc}
    \frac {\pm 2i \sin\pi\gamma(z)} {e^{\pm i\pi\gamma(z)}} &  \pm e^{- 2n  \phi_\beta (z)}     \\
    \mp e^{2n \phi_\beta (z)\mp 2i\pi \gamma(z)}& 1 \\
  \end{array}
\right),
   & \begin{array}{l}
   z\in \beta \pm i(0,\varepsilon).
     \end{array}
\end{array}
\right .
\end{equation}

\item[($T_c$)]      The asymptotic behavior of $T(z)$  at infinity
  is
  \begin{equation}\label{T-infty}T(z)= I+O\left (  1/ z\right ) ,~~~\mbox{as}~~z\rightarrow
                             \infty  .\end{equation}

\item[($T_d$)]
  $T(z)$ has the following behavior
  \begin{equation}\label{T-pm1}
 T(z)=O\left(\ln|z\mp 1|\right),\quad \mbox{as}~z\rightarrow \pm1.
  \end{equation}
\end{description}

To verify the jump conditions, we may use   the phase condition \eqref{g-phase-condition}, the relations \eqref{phi-alpha-g-relation} and \eqref{phi-beta-g-relation} between the $g$-function and the $\phi$-functions, and the following representations derived from them:
\begin{equation}\label{g-in-phi}
 g_+(x)- g_-(x)  =\left \{
 \begin{array}{ll}
 2\pi i, &x\in (-\infty, \alpha), \\[.2cm]
  2\pi i -2 \left (\phi_\alpha\right ) _+(x)= 2\pi i +2 \left (\phi_\alpha\right ) _-(x),& x\in (\alpha, 1),
  \\[.2cm]
  \frac {2\pi i} n \gamma(x)  -2 \left (\phi_\beta\right ) _+(x)= \frac {2\pi i} n \gamma(x)  +2 \left (\phi_\beta\right ) _-(x),& x\in (1, \beta),
   \\[.2cm]
  \frac {2\pi i} n \gamma(x),& x\in (\beta, \infty).
 \end{array}\right .
  \end{equation}

\subsection{The second transformation $T\to S$}

\begin{figure}[t]
\begin{center}
\includegraphics[width=12cm,bb=30 320 570 470]{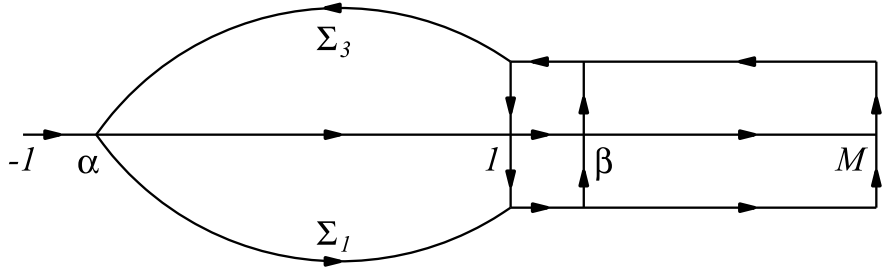}
 \caption{The contour $\Sigma_S$ for $S(z)$. We denote by $\Omega_+$ the  domain bounded by $\Sigma_3$ and the real axis such that $\alpha<\Re z<1$, and by  $\Omega_-$ the symmetric part bounded by $\Sigma_1$ and $(\alpha, 1)$.}
 \label{contour-for-S}
\end{center}
\end{figure}

Now we are in a position to apply the second transformation $T\to S$, associated with factorizations of the jump matrices.
        The transformation is defined explicitly as
\begin{equation}\label{transformation-T-to-S}
  S(z)=T(z)\left\{
  \begin{array}{ll}
   \left(
          \begin{array}{cc}
            1 &  0 \\
          \mp\phi_0(z) e^{2n \phi_\alpha (z)}   & 1\\
          \end{array}
        \right) & \mbox{for}~~z\in \Omega_\pm;~ \mbox{see~ Figure}~ \ref{contour-for-S}, \\[.6cm]
   \left(
          \begin{array}{cc}
            1 &  0 \\
          \pm e^{2n \phi_\beta (z)\mp 2\pi i \gamma(z)}   & 1\\
          \end{array}
        \right) & \mbox{for}~~z\in (1, \beta)\times (0, \pm i\varepsilon ), \\[.6cm]

              I & \mbox{otherwise}.
  \end{array}\right .
\end{equation}
We see that $S(z)$ solves the following RHP
\begin{description}
  \item[($S_a$)]     $S(z)$ is analytic in
  $\mathbb{C}\backslash \Sigma_S$, where the oriented contour $\Sigma_S$ is illustrated in Figure \ref{contour-for-S}.
  \item[($S_b$)] The jump conditions are
  \begin{equation}\label{S-jump}
  S_+(z)=S_-(z) J_S(z),
 ~~~ z\in \Sigma_S,\end{equation}with jumps on real axis
 \begin{equation}\label{S-jump-real-axis}
 J_S(x)=\left\{
  \begin{array}{ll}
    \left(
  \begin{array}{cc}
  \frac 1 { \left(\phi_0\right ) _-(x)}& e^{-2n\phi_\alpha(x)} \\
  0 & \frac 1 { \left(\phi_0\right ) _+(x)} \\
  \end{array}
\right),
   & x\in (-1, \alpha),
  \\[.6cm]
  \left(
    \begin{array}{cc}
      0 & 1 \\
      -1 & 0 \\
    \end{array}
  \right), &x\in (\alpha, 1)\cup(1, \beta),\\[.6cm]
 \left (\begin{array}{cc}
   1 & 0  \\
     -e^{2n\phi_\beta(x)}  &  1
 \end{array}\right ),& x\in (\beta, M),
   \end{array}\right .
   \end{equation}
and  jumps on the other contours
    \begin{equation}\label{S-jump-off}
 J_S(z)=\left\{
  \begin{array}{ll}
  \left(
    \begin{array}{cc}
      1 & 0 \\
      \phi_0(z) e^{2n\phi_\alpha(z)}  & 1 \\
    \end{array}
  \right), &z\in \Sigma_1,\\[.6cm]
 \left(
    \begin{array}{cc}
      1 & 0 \\
        -\phi_0(z) e^{2n\phi_\alpha(z)}  & 1 \\
    \end{array}
  \right), &z\in \Sigma_3,\\[.6cm]

  \left(
          \begin{array}{cc}
           1  & 0 \\
          \frac {e^{2n\phi_\beta (z)\mp \pi i \gamma(z)}} {2 i\sin \pi\gamma(z)}  & 1 \\
          \end{array}
        \right)
  ,  &  \left\{ \begin{array}{l}
                                                           z\in (\beta, M)\pm i\varepsilon, \\
                                                              z\in M\pm i(0,\varepsilon) ,
                                                          \end{array}\right .
  \\[.6cm]
\left (\begin{array}{cc}
   1 & \pm  e^{-2n\phi_\beta(z)}  \\
     \frac { e^{2n\phi_\beta(z)\mp \pi i\gamma(z)}}{2i \sin \pi\gamma(z)}  &  \frac { e^{ \pm \pi i\gamma(z)}}{\pm 2i \sin \pi\gamma(z)}
 \end{array}\right ),& z\in (1,\beta)\pm i\varepsilon,\\[.6cm]
   \left(
  \begin{array}{cc}
  1 & \pm  e^{-2n\phi_\beta(z)}\\
   0  &  1  \\
  \end{array}
\right),
   & z\in 1\pm i(0,\varepsilon) ,\\[.6cm]
 \left(
  \begin{array}{cc}
    1  & \pm e^{-2n\phi_\beta(z) }\\
   0 &  1\\
  \end{array}
\right),
   & z\in \beta \pm i(0,\varepsilon) .

 \end{array}\right .
   \end{equation}

\item[($S_c$)] The behavior at infinity is
 \begin{equation}\label{S-infty}
   S(z)=I+O(1/z), ~~z\to \infty.
 \end{equation}
\item[($S_d$)]
 $S(z)$ has the following behavior
  \begin{equation}
 S(z)=O\left(\ln|z\mp 1 |\right),\quad \mbox{as}~z\rightarrow \pm 1.
  \end{equation}
\end{description}

\subsection{The functions $\phi_\alpha$ and $\phi_\beta$}\label{phi-alpha-phi-beta}

Let us take a closer look at $\phi_\alpha(z)$ in a neighborhood of $z=\alpha$; cf. \eqref{phi-alpha}.

 {\lem{\label{phi-alpha-lemma}The function $\phi_\alpha(z)$ possesses the following convergent
  expansion
\begin{equation}\label{phi-alpha-expansion}
  n\phi_\alpha (z)=ic_\alpha(n)\tau_\alpha^{3/2}\sum^\infty_{k=0}c_{\alpha, k}(n) \tau_\alpha^k,~~ \tau_\alpha=\frac {z-\alpha}{\alpha+1},~~\arg \tau_\alpha\in (0, 2\pi),~~|\tau_\alpha|<1,
\end{equation}where $c_{\alpha, 0}=1$, and $c_\alpha(n)=\frac {\sqrt{\beta-\alpha}} 3 (1+\alpha)^{3/2}  A(\alpha)\sim \frac {2\sqrt{2b}} 3 \sqrt n$ for large $n$. It also holds the following asymptotic approximation for large $n$
 \begin{equation}\label{phi-alpha-in-n}
  n\phi_\alpha (z)= i\sqrt n \tau_\alpha^{3/2} \left (\xi_{0}(\tau_\alpha) +O\left (1/\sqrt n\right )\right )~~\mbox{for}~~|\tau_\alpha|=r,~~r\in (0,1),
\end{equation}
where $\xi_{0}(\tau)$ is an analytic function in $|\tau|<1$, depending only on $\tau$,  explicitly expressed in terms of the Beta function as
$\xi_{0}(\tau)=\sum^\infty_{k=0}   (-1)^k \sqrt{  b /2}\; B(k+2,    1 /2) \tau^k.$
  }}

 \vskip .2cm \noindent
{\bf{Proof: }}
 First, since $h_\alpha(z)$ is analytic in $\mathbb{C}\setminus (-\infty, -1]$, such that $h_\alpha(x)=-r'(x)/r(x)$ for $x\in (-1, 1)$; cf. \eqref{h-alpha},   from \eqref{G-solution} we can write
$$\pi i G(z)=\frac {\sqrt{(z-\alpha)(z-\beta)}}{2n}\left [h_\alpha(z) \int^\beta_\alpha \frac 1 {\pi\sqrt{(x-\alpha)(\beta-x)}}\frac {dx}{x-z}+A(z)\right ],$$  where $A(z)$, given in \eqref{A-def} below,  is real for $z\in (-1, \alpha)$,  and is analytic in $\mathbb{C}\setminus\{(-\infty, -1]\cup [1,\infty)\}$, to which  $\alpha$ belongs.
To justify the   analyticity of $A(z)$, one may use the equality
\begin{equation}\label{A Cauchy integral}
     \int^\beta_\alpha \frac 1 {\pi\sqrt{(x-\alpha)(\beta-x)}}\frac {dx}{x-z}=-\frac 1 { \sqrt{(z-\alpha)(z-\beta)}}~~\mbox{for}~z\not\in [\alpha, \beta],
\end{equation}
where branches are chosen such that $\arg(z-\alpha)\in (-\pi, \pi)$ and  $\arg(z-\beta)\in (-\pi, \pi)$.

Substituting all above into \eqref{nu-representation}, we have
\begin{equation}\label{phi-alpha-derivative}
 \phi_\alpha'(z)=\nu_\alpha(z)=\frac {\sqrt{(z-\alpha)(z-\beta)}} {2n} A(z)=   i \sqrt{z-\alpha} A_\alpha(z)
\end{equation}
 in a neighborhood of $z=\alpha$, with  $\arg (z-\alpha)\in (0, 2\pi)$, and $A_\alpha(z)$ is  analytic at $z=\alpha$, being real for $z\in (-1, \alpha)$. A similar discussion can be found in Zhou and Zhao \cite{zhou-zhao}.

Hence, we accordingly derive \eqref{phi-alpha-expansion}  from \eqref{phi-alpha}.
We include some details for the above approximation. Recalling that
\begin{equation}\label{A-def}
 A(z)= \int^\beta_\alpha \frac  {-h_\alpha(z)   -\frac {r'(x)}{r(x)}} {x-z}  \frac {dx}  {\pi\sqrt{(x-\alpha)(\beta-x)}}+\int^\infty_\beta \frac {2bx(x^2-1)^{-3/2}}{ \sqrt{(x-\alpha)(x-\beta)}} \frac {dx}{x-z}
\end{equation}
for $z\in \mathbb{C}\setminus\{(-\infty, -1]\cup [1,\infty)\}$; cf. \eqref{G-solution}.
By careful estimation similar to that in Section \ref{MRS}, for large $n$, we have
$$A(\alpha) \sim\int^0_\alpha \frac  {h_\alpha(x)-h_\alpha(\alpha)}{x-\alpha}  \frac {dx}  {\pi\sqrt{(x-\alpha)(\beta-x)}}
\sim \frac 1 {\sqrt2 \pi}\int^0_\alpha    \left [h_\alpha(x)-h_\alpha(\alpha)\right  ]  \frac {dx} {(x-\alpha)^{3/2}},$$
with an error
$O\left(n^{3/2}\right )$. Here
use has been made of \eqref{r-end-behavior} and the fact that $\alpha\sim -1+b/n$ for large $n$. Now integrating by parts once, we further have
$$A(\alpha)\sim   \frac  {\sqrt 2}{ \pi}\int^0_\alpha      \frac {h_\alpha'(x) dx} {\sqrt {x-\alpha}}\sim \frac {3b} 2 \int^0_\alpha \frac {dx}{(1+x)^{5/2}\sqrt{x-\alpha}}\sim 2b(1+\alpha)^{-2}, $$
each time with an error  $O\left(n^{3/2}\right )$. Here we have taken into account the fact that
$$h_\alpha'(x)=  \frac {3\pi b}{2\sqrt 2} (1+x)^{-5/2}+ O\left ((1+x)^{-2}\right )~~\mbox{as}~~x\to -1^+; $$ cf. \eqref{h-alpha}  and \eqref{r-end-behavior}.   Hence we obtain the above approximation for  $c_\alpha(n)$.

Now we turn to the evaluation of  the leading behavior of $n\phi_\alpha (z)$ for large $n$ and  mild $\tau=\frac {z-\alpha}{1+\alpha}$. To this aim, we may go further to approximate $A^{(k)}(\alpha)$ for $\alpha=1,2,\cdots.$ Similar to the approximating of $A(\alpha)$, we have
$$A^{(k)}(\alpha)\sim \frac {k!} {\sqrt 2\; \pi} \int^0_\alpha \frac {h_\alpha(x)-h_\alpha(\alpha)-h_\alpha'(\alpha) (x-\alpha) -\cdots - \frac {h_\alpha^{(k)}(\alpha)} {k!} (x-\alpha)^k} {(x-\alpha)^{k+\frac 3 2}} dx.$$
Integrating by parts $k+1$ times, we have
$$\frac {A^{(k)}(\alpha)}{k!} \sim \frac {\Gamma(\frac 1 2)} {\sqrt 2\; \pi \Gamma(k+\frac 3 2) }\int^0_\alpha \frac {h_\alpha^{(k+1)}(x) dx}{\sqrt{x-\alpha}}\sim b  (-1)^k \left (k+\frac 3 2\right ) \int^0_\alpha \frac {dx}{(1+x)^{k+\frac 52} \sqrt{x-\alpha}}.$$
Here use has been made of the fact that
$$h_\alpha^{(k+1)}(x)\sim (-1)^k   \frac {\pi b}{ \sqrt 2} \frac 3 2 \frac 5 2\cdots\left  (k+\frac 3 2\right ) (1+x)^{-k-\frac 5 2}~~\mbox{as}~~x\to -1^+; $$
cf. \eqref{h-alpha}. Extending the integral interval to $[\alpha, \infty)$, one obtains
$$\frac {A^{(k)}(\alpha)}{k!} \sim (-1)^{k}   \frac {b \left (k+\frac 3 2\right )  B(k+2,  \frac 1 2)}{(1+\alpha)^{k+2}}~~\mbox{for}~~k=0,1,2,\cdots .$$

Expanding $A(z)$ into a Maclaurin series in $\tau=\frac {z-\alpha}{1+\alpha}$,   substituting all above to \eqref{phi-alpha-derivative} and integrating it, we obtain \eqref{phi-alpha-in-n}, thus completing the proof.\hfill\qed\vskip.4cm

Now we see from \eqref{phi-alpha-expansion} that in a neighborhood of $z=\alpha$ cutting alone $[\alpha, \infty)$, $e^{2n\phi_\alpha(z)}$ is exponentially small for $\arg (z-\alpha)\in (0, 2\pi/3)\cup (4\pi/3, 2\pi)$, and is exponentially large for $\arg (z-\alpha)\in (2\pi/3, 4\pi/3)$, so long as $|\tau|=\left | \frac {z-\alpha}{\alpha+1}\right |$ lies in compact subsets of $(0,1)$.

We proceed to show that the jumps $J_S$ for $S$ off the real line are of the form $I$ plus exponentially small terms; cf \eqref{S-jump-off}. This may be roughly explained by using the Cauchy-Riemann condition. For example, in view of \eqref{phi-alpha} and \eqref{phi-beta}, we have negative derivatives
$$\frac {\partial \Re \phi_\alpha}{\partial y}=-\pi \psi(x),~~\frac {\partial \Re \left\{\phi_\beta-\frac{\pi i\gamma} n \right \}}{\partial y}=-\pi \psi(x),~~\mbox{and}~~-\frac {\partial \Re \phi_\beta}{\partial y}=-\pi\left [ \frac {-\gamma'(x)} n- \psi(x)\right ],$$respectively on the upper edge of $x\in (\alpha, 1)$, $x\in (1, \infty)$, and $x\in (1, \beta)$.

Since the contours depend on $n$, it is necessary to estimate the exponentials carefully.  An approach is to check the derivatives such as $\frac {\partial \Re \phi_\alpha}{\partial y}$ in neighborhoods of real segments, as carried out  in  \cite{zhou-zhao} and \cite{zhou-xu-zhao}. To show that the jumps \eqref{S-jump-off} off the real line are of the form $I$ plus exponentially small terms, we use an alternative straightforward way here.    Indeed, substituting \eqref{g-integral}  into \eqref{phi-alpha-g-relation}, we obtain
\begin{align}\label{phi-alpha-integral}
\phi_\alpha(z) &=\frac {\sqrt{(z-\alpha)(z-\beta)}} n \int_\beta^\infty\frac{-\gamma(x)}{ \sqrt{(x-\alpha)(x-\beta)}}\frac{dx}{x-z}- \ln\frac{ z-\frac {\alpha+\beta} 2 +\sqrt{(z-\alpha)(z-\beta)}}{\frac {\beta-\alpha} 2}\nonumber\\
&~~ -\frac{\sqrt{(z-\alpha)(z-\beta)}}{2n\pi}\int_\alpha^\beta\frac{\ln \pi r(x) }{\sqrt{(x-\alpha)(\beta-x)}}\frac{dx}{x-z}-\frac {ib\ln(z+i\sqrt{1-z^2})}{n\sqrt{1-z^2}}\pm \pi i.
\end{align}We take  part of $\Sigma_3$ (cf. Figure \ref{contour-for-S}) as an example to show $e^{2n\phi_\alpha(z)}$ is exponentially small there.  Along  the slope segment  $\arg (z-\alpha)=\theta_+$ and $\alpha+r(1+\alpha)\leq \Re z\leq \alpha +\frac {\delta_+}{\tan\theta_+}$, we may pick the dominant  contribution from the logarithm in \eqref{phi-alpha-integral}, so that $\Re\{2n\phi_\alpha(z)\}\leq-  C n \sqrt{|z-\alpha|}$. Other exponential terms can be estimated similarly using the explicit representation.

Analysis can also be carried out for $\phi_\beta(z)$ in an $O(1/n)$ neighborhood of $z=\beta$. A new variable used here is $\tau_\beta=\frac {z-\beta}{\beta-1}$. Indeed, a combination of \eqref{G-solution} and \eqref{A Cauchy integral} with \eqref{nu-beta} gives
\begin{equation}\label{phi-beta-derivative}
  \phi_\beta'(z)=\nu_\beta(z)=\frac {\sqrt{(z-\alpha)(z-\beta)}}{2n}\; B(z),
\end{equation}where
\begin{equation}\label{B-def}
   B(z)= \int^\beta_\alpha \frac  {\frac {r'(z)}{r(z)} -\frac {r'(x)}{r(x)}} {x-z}  \frac {dx}  {\pi\sqrt{(x-\alpha)(\beta-x)}}-\int^\infty_\beta \frac {2\gamma'(x)}{ \sqrt{(x-\alpha)(x-\beta)}} \frac {dx}{x-z}\pm
   \frac {2\pi i \gamma'(z)}{\sqrt{(z-\alpha)(z-\beta)}}
\end{equation} for $\pm\Im z>0$, with $r(x)$ given in \eqref{r-def} for real $x$, $r(z)$ is the function appeared  in \eqref{nu-beta}, analytic  in  $\mathbb{C}\setminus (-\infty, 1]$, and $\gamma(z)=\frac b{\sqrt{z^2-1}}$, analytic in $\mathbb{C}\setminus [-1, 1]$. It is readily verified that $B(z)$ is analytic in  $\mathbb{C}\setminus (-\infty, 1]$.

For large $n$, the leading order contribution to $B(\beta)$ comes from the last two terms in \eqref{B-def}. Actually we have
$$B(\beta)\sim  \sqrt 2 \int^\infty_\beta \frac {\gamma'(\beta)-\gamma'(x)}{(x-\beta)^{3/2}} dx= \int^\infty_\beta \frac {-2^{\frac 3 2}\gamma''(x) dx}{\sqrt{x-\beta}}\sim \int^\infty_\beta \frac {-\frac {3b} 2 dx}{(x-1)^{5/2} \sqrt{x-\beta}}=\frac {-2b}{(\beta-1)^{2}}.$$
Hence in analog  to Lemma \ref{phi-alpha-lemma}, we have
{\lem{\label{phi-beta-lemma}
The function $\phi_\beta(z)$ possesses the following convergent
  expansion
\begin{equation}\label{phi-beta-expansion}
  n\phi_\beta (z)=-c_\beta(n)\tau_\beta^{3/2}\sum^\infty_{k=0}c_{\beta, k}(n) \tau_\beta^k,~~\tau_\beta=\frac {z-\beta}{\beta-1},~~ \arg\tau_\beta\in (-\pi, \pi),~~|\tau_\beta|<1,
\end{equation}where   $c_{\beta, 0}=1$, and $c_\beta(n) \sim \frac {2\sqrt{2b}} 3 \sqrt n$ for large $n$. Also, similar to  \eqref{phi-alpha-in-n}, it holds
\begin{equation}\label{phi-beta-in-n}
  n\phi_\beta (z)= -\sqrt n \tau_\beta^{3/2} \left (\xi_{0}(\tau_\beta) +O\left (1/\sqrt n\right )\right )~~\mbox{for}~~|\tau_\beta|=r,~~r\in (0,1),~~n\to\infty,
\end{equation}where $\xi_{0}(\tau)$ is the same function as in \eqref{phi-alpha-in-n},  analytic   in the unit disc.\hfill \qed}}

\subsection{Parametrix for the outside region} \label{outside-region}

For fixed $z\in \Sigma_S$, with the possible exceptions of $(-1, \alpha)$ and $(\alpha, \beta)$, it can be verified that $J_S(z)$ is $I$, plus an exponentially small term, for  large-$n$.   Hence, we have the limiting RH problem:

\begin{description}
  \item[($N_a$)]     $N(z)$ is analytic in
  $\mathbb{C}\backslash [-1, \beta]$.
  \item[($N_b$)] The jump condition is
  \begin{equation}\label{N-jump}
  N_+(x)=N_-(x) \left \{
  \begin{array}{ll}
\left(
  \begin{array}{cc}
   { \left(\phi_0\right ) _+(x)}& 0 \\
  0 & \frac 1 { \left(\phi_0\right ) _+(x)} \\
  \end{array}
\right), &  x\in (-1, \alpha),  \\[.5cm]
    \left(
    \begin{array}{cc}
      0 & 1 \\
      -1 & 0 \\
    \end{array}
  \right), &  x\in (\alpha, \beta),
\end{array}\right .
  \end{equation}
  where $\phi_0(z)$ is defined in \eqref{phi-0}, such that $\left(\phi_0\right ) _+(x)\left(\phi_0\right ) _-(x)=1$ for $x\in (-1, \alpha)$.

 \item[($N_c$)] The behavior at infinity is
 \begin{equation}\label{N-infty}
   N(z)=I+O(1/z), ~~z\to \infty.
 \end{equation}

\end{description}

A solution can be constructed
explicitly as
\begin{equation}\label{N-expression}
 N(z) =D(\infty)^{-\sigma_3} \left(
         \begin{array}{cc}
           1+  \frac{f_s( \tau_\beta )} 2 &-\frac {i  f_s(  \tau_\beta)} 2 \\
          -\frac{i   f_s( \tau_\beta)} 2 &1- \frac { f_s(  \tau_\beta)} 2 \\
         \end{array}
       \right)  \left(
         \begin{array}{cc}
           1-  \frac{f_s( \tau_\alpha )} 2 & -\frac {i  f_s(  \tau_\alpha)} 2 \\
          -\frac{i   f_s( \tau_\alpha)} 2 &1+ \frac { f_s(  \tau_\alpha)} 2 \\
         \end{array}
       \right) N_0(z) D(z)^{\sigma_3},
 \end{equation}where $f_s( \tau)=\frac 1 {12 b} \frac 1 \tau+ \frac 5 {48b} \frac 1 {\tau^2}$, $\tau_\alpha=\frac {z- \alpha}{1+\alpha}$, $\tau_\beta=\frac {z- \beta}{ \beta-1}$, $D(z)$ is the Szeg\"{o} function, namely, a function analytic and non-vanishing  in $\{\mathbb{C}\setminus [-1, \beta]\}\cup\{\infty\}$, such that $D_+(x)D_-(x)=1$ for $x\in (\alpha, \beta)$ and  $D_+(x)/D_-(x)=\left (\phi_0\right )_+(x)$ for $x\in (-1, \alpha)$; see \cite{KuiVanAssche}  for a discussion of the  Szeg\"{o} functions, and \cite{xdz2014} for a relevant  construction. In the present  case, we have
$$D(z)=\exp\left ({-\frac {\sqrt{(z-\alpha)(z-\beta)}}{2\pi i} \int^\alpha_{-1} \frac {\ln\left ( \left(\phi_0\right ) _+(x)\right )}{\sqrt{(\alpha-x)(\beta-x)}}\frac {dx}{x-z}}\right ),$$
and
 $$    N_0(z)=
\left(
  \begin{array}{cc}
    \frac {\varrho(z) + \varrho^{-1}(z)} {2}&\frac {\varrho(z) - \varrho^{-1}(z)} {2i} \\
    -\frac {\varrho(z) - \varrho^{-1}(z)} {2i} &\frac {\varrho(z) + \varrho^{-1}(z)} {2} \\
  \end{array}
\right)  ~~\mbox{with}~~
 \varrho(z)=\left ( \frac {z-\beta}{z-\alpha} \right )^{\frac 1 4},$$  where
$\arg(z-\alpha)\in (-\pi, \pi)$ and $\arg(z-\beta)\in (-\pi, \pi)$. Thus for large $n$,
$$ D(\infty)=
\exp\left ({ \frac {1}{2\pi i} \int^\alpha_{-1} \frac {\ln\left ( \left(\phi_0\right ) _+(x)\right )dx}{\sqrt{(\alpha-x)(\beta-x)}}}\right )\sim  \exp\left (\frac {\left (1-\frac 1 {6b}\right )}{2\pi}
\int^\alpha_{-1} \frac {\sqrt {x+1} \; dx}{\sqrt{ \alpha-x}}\right )=1+O\left(\frac 1 n\right ). $$

It is worth noting that the  extra  rational factors on the left are needed to accomplish the matching conditions below for $P^{(\alpha)}$ and $P^{(\beta)}$. (see \cite[(4.15)]{zhou-zhao}, and \cite{zhou-xu-zhao}).

\subsection{Local parametrix at $z=\alpha$}

The limiting RHP $N(z)$ fails to approximate $S(z)$ at $z=\alpha$ and $z=\beta$ since the jump  $J_S(z)$ is not close to $I$ as $z\to \alpha$. We need to construct the following parametrix at a neighborhood $|\tau|<r$, denoted by $U_r(\alpha)$,  of   $z=\alpha$. Here we have used a re-scaled variable $\tau=\frac {z-\alpha}{1+\alpha}$.

\begin{figure}[t]
\begin{center}
\includegraphics[width=6cm]{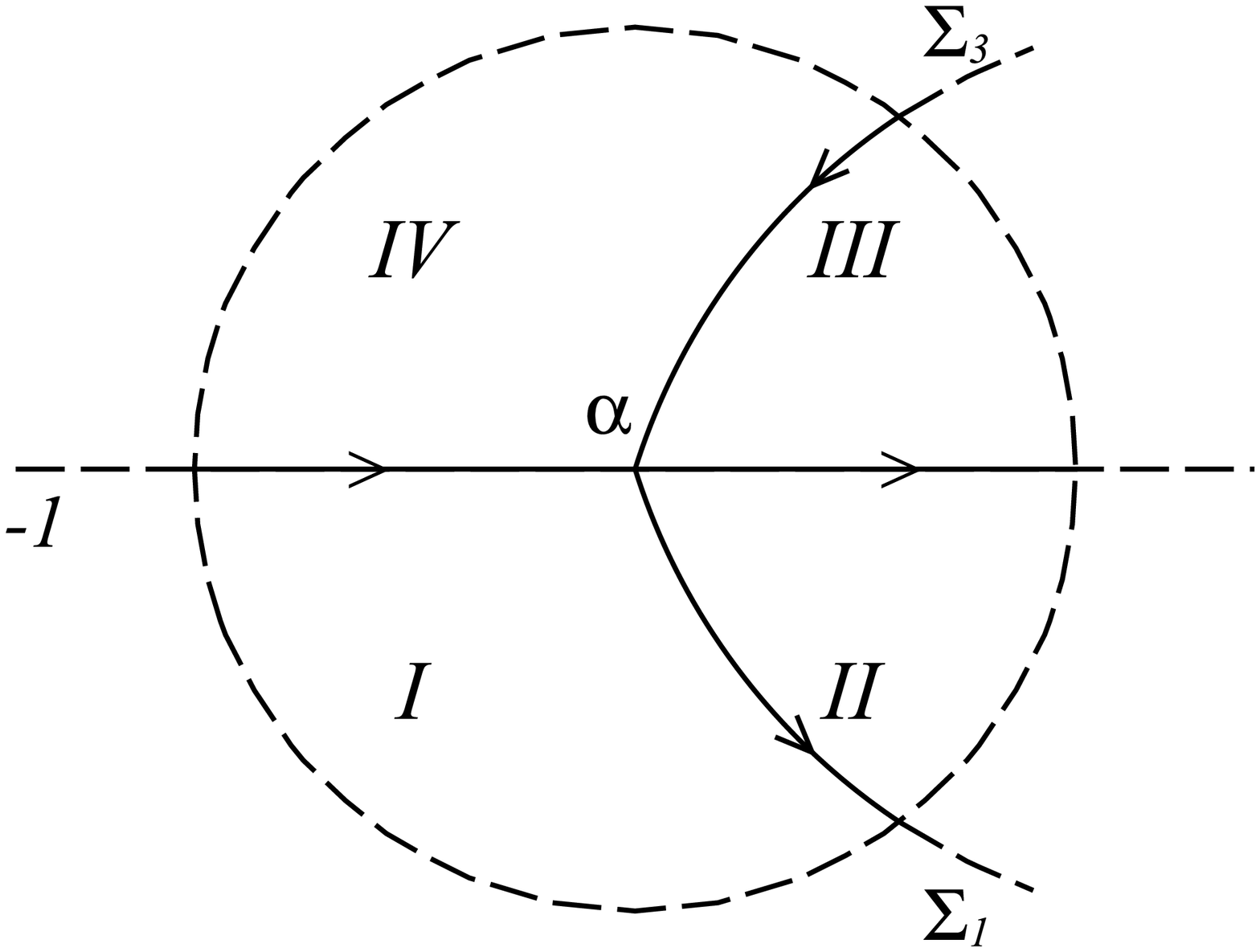}\hskip 1cm \includegraphics[width=6cm]{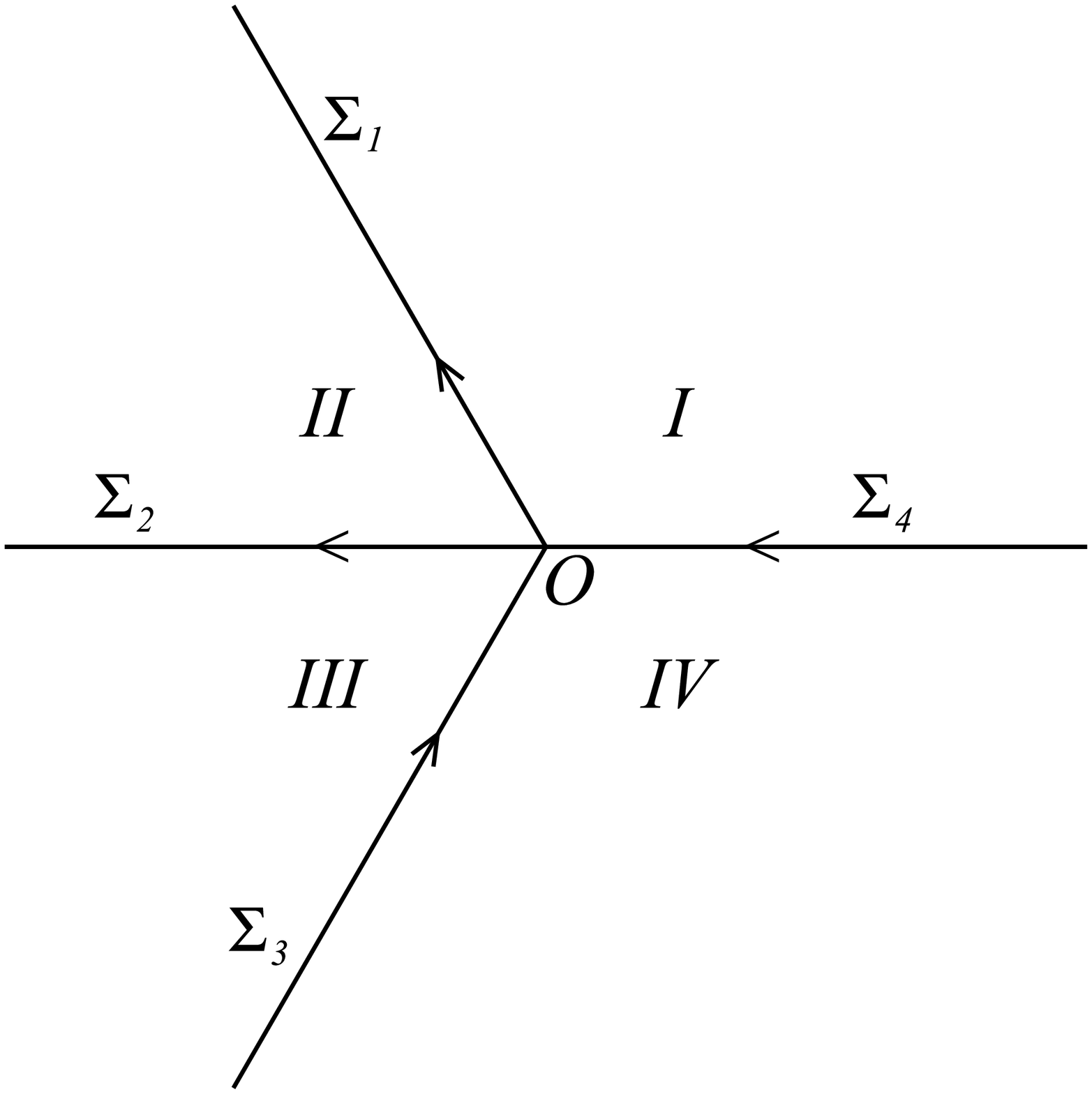}
 \caption{The contours for $P^{(\alpha)}(z)$ (left), and the contours and sectors for the model RHP   $\Psi(s)$ (right), connected by the conformal mapping $s=\lambda_\alpha(z)$. We still denote by   $\Sigma_1$ and $\Sigma_3$ the $s$-images (right) of the   oriented $z$-contours (left).  }
 \label{contour-for-P-alpha}
\end{center}
\end{figure}

 \begin{description}
  \item[($P^{(\alpha)}_a$)]     $P^{(\alpha)}(z)$ is analytic in $U_r(\alpha)\setminus \Sigma_S$; cf. Figures \ref{contour-for-S} and \ref{contour-for-P-alpha} for $\Sigma_S$.
  \item[($P^{(\alpha)}_b$)]  $P^{(\alpha)}(z)$ satisfies the same  jump conditions as $S(z)$ in  $U_r(\alpha)$, namely,
  \begin{equation}\label{P-alpha-jump}
  P^{(\alpha)}_+(z)=P^{(\alpha)}_-(z)
  \left\{
  \begin{array}{ll}
  \left(
  \begin{array}{cc}
  \frac 1 {\left (\phi_0\right )_-(z)} & e^{-2n\phi_\alpha(z)} \\
  0 &  \frac 1 {\left (\phi_0\right )_+(z)}  \\
  \end{array}
\right),
   & z\in (-1, \alpha)\cap U_r(\alpha),
  \\[.6cm]
  \left(
    \begin{array}{cc}
      0 & 1 \\
      -1 & 0 \\
    \end{array}
  \right), &z\in (\alpha, 1)\cap U_r(\alpha),\\[.6cm]
 \left(
    \begin{array}{cc}
      1 & 0 \\
      \phi_0(z) e^{2n\phi_\alpha(z)}  & 1 \\
    \end{array}
  \right), &z\in \Sigma_1\cap U_r(\alpha), \\ [.6cm]
   \left(
    \begin{array}{cc}
      1 & 0 \\
        -  \phi_0(z) e^{2n\phi_\alpha(z)}  & 1 \\
    \end{array}
  \right), &z\in \Sigma_3\cap U_r(\alpha).
     \end{array}\right .
   \end{equation}

 \item[($P^{(\alpha)}_c$)] The matching condition holds
 \begin{equation}\label{P-alpha-matching}
  P^{(\alpha)}(z)=\left ( I+O\left (\frac 1 n\right ) \right ) N(z), ~~z\in \partial U_r(\alpha).
 \end{equation}
\end{description}

A solution to the RHP has been constructed in, e.g., \cite{dkmv2}, see also \cite{zhou-zhao}, expressed in terms of the following matrix-valued function
\begin{equation}\label{Airy}
  \Psi(s)=\left\{
                 \begin{array}{ll}
                  \left(
                     \begin{array}{cc}
                        \Ai(s) & \Ai(\omega^2 s) \\
                        \Ai'(s)&\omega^2 \Ai'(\omega^2 s) \\
                      \end{array}
                  \right)e^{-\frac{\pi i}{6}\sigma_3},  &  s\in I ;  \\[.5cm]
                       \left(
                     \begin{array}{cc}
                        \Ai(s) & \Ai(\omega^2 s) \\
                        \Ai'(s)&\omega^2 \Ai'(\omega^2 s) \\
                      \end{array}
                  \right)e^{-\frac{\pi i}{6}\sigma_3}\left(
                                                       \begin{array}{cc}
                                                         1 & 0 \\
                                                         -1 & 1 \\
                                                       \end{array}
                                                     \right),
                    &  s\in II ;  \\[.5cm]
                            \left(
                     \begin{array}{cc}
                        \Ai(s) & -\omega^2 \Ai(\omega s) \\
                        \Ai'(s)&- \Ai'(\omega s) \\
                      \end{array}
                  \right)e^{-\frac{\pi i}{6}\sigma_3}\left(
                                                       \begin{array}{cc}
                                                         1 & 0 \\
                                                      1 & 1 \\
                                                       \end{array}
                                                     \right),
                    &  s\in III ;  \\[.5cm]
                        \left(
                     \begin{array}{cc}
                        \Ai(s) & -\omega^2 \Ai(\omega s) \\
                        \Ai'(s)&- \Ai'(\omega s) \\
                      \end{array}
                  \right)e^{-\frac{\pi i}{6}\sigma_3},  &  s\in IV ;
                 \end{array}
          \right.
 \end{equation}
 cf. \cite[(7.9)]{dkmv2}, where  $\omega=e^\frac{2\pi i}{3}$, and the sectors $I-IV$ are illustrated in Figure \ref{contour-for-P-alpha}. We note that with the orientation  of $\Sigma_k$ indicated, $\Psi(s)$ possesses the constant jumps
$$J_\Psi(s)=\left\{
  \begin{array}{ll}
  \left(
  \begin{array}{cc}
  1 & 0 \\
  -1 & 1 \\
  \end{array}
\right),
   & s\in \Sigma_1,
  \\[.4cm]
  \left(
    \begin{array}{cc}
      0 & -1 \\
      1 & 0 \\
    \end{array}
  \right), &s\in \Sigma_2,\\[.4cm]
 \left(
    \begin{array}{cc}
      1 & 0 \\
       1   & 1 \\
    \end{array}
  \right), &s\in \Sigma_3, \\ [.4cm]
   \left(
    \begin{array}{cc}
      1 &-1 \\
       0  & 1 \\
    \end{array}
  \right), &s\in \Sigma_4.
     \end{array}\right .
$$ To verify the jump conditions, use has been made of the fact that
\begin{equation}\label{Airy-connection}
 \Ai(s)+\omega\Ai(\omega s)+\omega^2\Ai(\omega^2 s)=0~~\mbox{for}~~s\in \mathbb{C}.
\end{equation}

Resuming the re-scaled variable $\tau=\tau_\alpha=\frac {z-\alpha}{1+\alpha}$, we see from \eqref{phi-alpha-expansion} that
\begin{equation}\label{alpha-conformal-mapping}
  \lambda_\alpha(z)=  e^{-4\pi i/3}  \left ( \frac 3 2\right )^{2/3} n^{1/3} \phi^{2/3}_\alpha (z)
\end{equation}defines a conformal mapping from an $O(1)$ neighborhood of $\tau_\alpha=0$ to  an $O(1)$ neighborhood of $\lambda_\alpha=0$.
We can then write down the parametrix as
\begin{equation}\label{P-alpha-representation}
  P^{(\alpha)}(z)=E(z) \Psi\left (n^{1/3}\lambda_\alpha(z)\right ) e^{-\frac \pi 2 i\sigma_3}e^{n\phi_\alpha(z)\sigma_3}\left (\phi_0(z)\right )^{\frac 1 2\sigma_3}, ~~z\in U_r(\alpha),
\end{equation}where $\phi_0^{\frac 1 2}$ is determined such that $ ( \phi_0^{\frac 1 2} )_\pm (x) =e^{\pm \frac 1 2 i\Theta(x)}$, with  $\left ( \phi_0 \right )_\pm (x) =e^{\pm  i\Theta(x)}$ for $x\in (-1, 1)$ and $\Theta(x)$ being real,  $E(z)$ is an analytic function in $z\in U_r(\alpha)$, as can be determined by the matching condition \eqref{P-alpha-matching}, and by expanding $\Psi$ in \eqref{P-alpha-representation} for large $n$.

Indeed, as $\tau_\alpha=\frac {z-\alpha}{1+\alpha}\sim O(1)$, $s:=n^{1/3} \lambda_\alpha(z) \sim O(n^{1/3})$ large, and $\zeta=\frac 2 3s^{3/2}=n\phi_\alpha(z) \sim e^{\frac {\pi i} 2} \sqrt n \tau_\alpha^{3/2} \xi_0(\tau_\alpha)$ for large $n$. Substituting in the   asymptotic approximation for $\Psi(s)$ in \eqref{Airy} gives
\begin{equation}\label{Psi-expansion}
     \Psi\left (s \right ) e^{-\frac \pi 2 i\sigma_3}e^{n\phi_\alpha(z)\sigma_3}\sim \frac 1 {2\sqrt \pi}
e^{-\frac \pi 6 i} {s^{-\frac 1 4\sigma_3}}  \left\{ \sum_{k=0}^\infty \frac {\Psi_{k}}{\zeta^k}\right\} M_\alpha;
\end{equation}
 see \cite[(10.4.7), (10.4.59), (10.4.61)]{as}, see also \cite[Sec.\;5.1]{zhou-xu-zhao} for a detailed discussion.
 In \eqref{Psi-expansion}, the coefficients are given as
 $$M_\alpha=\left(
       \begin{array}{cc}
         -i & -1 \\
        i & -1 \\
       \end{array}
     \right),~~\Psi_{2k}=\left(
                                   \begin{array}{cc}
                                     c_{2k} & 0 \\
                                     0 & d_{2k} \\
                                   \end{array}
                                 \right),~~\Psi_{2k+1}=\left(
                                   \begin{array}{cc}
                                    0 & c_{2k+1}  \\
                                   d_{2k+1} &  0  \\
                                   \end{array}
                                 \right)~~\mbox{for}~~k=0,1,2,\cdots,$$
with $c_0=d_0=1$, $c_k=\frac {\Gamma(3k+\frac 1 2)}{(54)^k k! \Gamma(k+\frac 1 2)}$ and $d_k=-\frac {6k+1}{6k-1} c_k$ for $k=1,2,\cdots$.
which shares the same jump condition as $N(z)$. Also, we need to go further to analyze the piece-wise analytic function  $\phi_0(z)$. From \eqref{phi-0}, and using Stirling's formula (see \cite[(6.1.37)]{as}), we
 obtain
 $$\phi_0(z)=1\pm \frac {ib_\alpha\sqrt{ \tau+1}}{\sqrt n} +O\left ( \frac 1 n\right ),~~\pm \Im z>0$$
    for $\tau_\alpha=\frac {z-\alpha}{\alpha+1}=O(1)$ and $n$ large, where $b_\alpha=\sqrt{2b}\left(1-\frac 1 {6b}\right )$, and $\arg (\tau+1)\in (-\pi, \pi)$.

For $\tau_\alpha=\frac{z-\alpha}{1+\alpha}=O(1)$,   we can rewrite $N(z)$ as
$$N(z)=\left (I+O\left (\frac 1 n\right )\right )  M_\alpha^{-1} \left(
                       \begin{array}{cc}
                         1 & f_s(\tau_\alpha) \\
                         0 & 1 \\
                       \end{array}
                     \right)M_\alpha   M_\alpha^{-1} \varrho(z)^{\sigma_3} M_\alpha D(z)^{\sigma_3}.$$
Then, we  can determine the analytic function in $z\in U_r(\alpha)$ as
\begin{equation}\label{E}
  E(z)=2\sqrt\pi e^{\frac \pi 6 i} M_\alpha^{-1}  \left(
                                                    \begin{array}{cc}
                                                      1 &  - 1+\frac 1{6b}- f_r(\tau_\alpha) \\
                                                      0 & 1 \\
                                                    \end{array}
                                                  \right) \left (s^{\frac 1 4} \varrho(z)\right )^{\sigma_3}
  ,
\end{equation}where $f_r(\tau)=f(\tau)-f_s(\tau)$ and $f_s(\tau)=\frac 1 {12 b} \frac 1 \tau+ \frac 5 {48b} \frac 1 {\tau^2}$ is   the  singular part of      $f(\tau)=\frac {c_1 \sqrt{2/b}}{\tau^2\xi_0(\tau)}$; cf. \eqref{phi-alpha-in-n}.
The matching condition \eqref{P-alpha-matching} is readily verified.
In the verification, use may be made of the facts that $M_\alpha \sigma_3 M_\alpha^{-1} =-\sigma_1$, and that for $|\tau|=O(1)$, $\ln \phi_0(z)=\pm 2c\sqrt{z+1}+O(1/n)$ respectively as $\pm \Im z>0$, with $c=\frac {ib_\alpha}{2\sqrt b}=\frac i {\sqrt 2}\left (1-\frac 1 {6b}\right )$. Hence
$$\phi_0^{\frac 1 2 }D(z)^{-1}
\sim \exp\left ({ \pm c  \sqrt{z+1}\pm \frac {c\sqrt{z-\alpha}}{ \pi} \int^\alpha_{-1} \frac {\sqrt{x+1}}{\sqrt{\alpha-x}}\frac {dx}{x-z}}
\right )=e^{  \pm c \sqrt{z-\alpha}} $$ for  $\pm \Im z>0$,  up to  an error of order $O(1/n)$,
as follows from the identity
$$-\frac 1\pi \int^\alpha_{-1} \frac {\sqrt{x+1}}{\sqrt{\alpha-x}}\frac {dx}{x-z}=\sqrt{\frac {z+1}{z-\alpha}}-1,~~z\in \mathbb{C},~~\arg(z+1), ~\arg(z-\alpha)\in (-\pi, \pi).$$

\subsection{Local parametrix at $z=\beta$}

We proceed to construct a  parametrix at a neighborhood $U_r(\beta)$ of   $z=\beta$, described  as $|\tau|<r$ for  $\tau=\frac {z-\beta}{\beta-1}$. We note that $\beta$ is the   band-saturated region endpoint.
The parametrix is formulated as
 \begin{description}
  \item[($P^{(\beta)}_a$)]     $P^{(\beta)}(z)$ is analytic in $U_r(\beta)\setminus \Sigma_S$; cf. Figure  \ref{contour-for-S}  for $\Sigma_S$.
  \item[($P^{(\beta)}_b$)]  $P^{(\beta)}(z)$ satisfies the same  jump conditions as $S(z)$ in  $U_r(\beta)$, namely,
  \begin{equation}\label{P-alpha-jump}
  P^{(\beta)}_+(z)=P^{(\beta)}_-(z)
  \left\{
  \begin{array}{ll}
  \left(
  \begin{array}{cc}
  1 & 0 \\
  - e^{2n\phi_\beta(z)} & 1 \\
  \end{array}
\right),
   & z\in (\beta, M)\cap U_r(\beta),
  \\[.6cm]
  \left(
    \begin{array}{cc}
      0 & 1 \\
      -1 & 0 \\
    \end{array}
  \right), &z\in (1, \beta)\cap U_r(\beta),\\[.6cm]
 \left(
    \begin{array}{cc}
      1 &  e^{-2n\phi_\beta(z)}\\
     0 & 1 \\
    \end{array}
  \right), &z\in \{\beta+i(0,\varepsilon)\}\cap U_r(\beta), \\ [.6cm]
   \left(
    \begin{array}{cc}
      1 & -e^{-2n\phi_\beta(z)} \\
       0 & 1 \\
    \end{array}
  \right), &z\in \{\beta-i(0,\varepsilon)\}\cap U_r(\beta).
     \end{array}\right .
   \end{equation}

 \item[($P^{(\beta)}_c$)] The matching condition holds
 \begin{equation}\label{P-beta-matching}
  P^{(\beta)}(z)=\left ( I+O\left (\frac 1 n\right ) \right ) N(z), ~~z\in \partial U_r(\beta).
 \end{equation}
\end{description}
The RHP has also been solved in earlier literature. Indeed, as in   Bleher and  Liechty \cite[(11.6)-(11.7)]{Bleher}, we use the matrix function
\begin{equation}\label{Psi-beta}
  \Psi_\beta(s)=\left\{
  \begin{array}{ll}
    \left(
       \begin{array}{cc}
         \omega^2\Ai(\omega^2 s) & -\Ai(s) \\
         \omega \Ai'(\omega^2 s) & -\Ai'(s)\\
       \end{array}
     \right)
     & \mbox{for}~\arg s\in (0,\pi/2), \\[.4cm]
    \left(
       \begin{array}{cc}
         \omega^2\Ai(\omega^2 s) & \omega \Ai(\omega s) \\
         \omega \Ai'(\omega^2 s) & \omega^2 \Ai'(\omega s)\\
       \end{array}
     \right)
     & \mbox{for}~\arg s\in (\pi/2, \pi), \\[.4cm]
       \left(
       \begin{array}{cc}
        \omega \Ai(\omega s) &  -\omega^2\Ai(\omega^2 s)  \\
         \omega^2 \Ai'(\omega s) &- \omega \Ai'(\omega^2 s)\\
       \end{array}
     \right)
     & \mbox{for}~\arg s\in (-\pi,-\pi/2), \\[.4cm]
    \left(
       \begin{array}{cc}
         \omega\Ai(\omega s) & \Ai(  s) \\
         \omega^2 \Ai'(\omega s) &   \Ai'(  s)\\
       \end{array}
     \right)
     & \mbox{for}~\arg s\in (-\pi/2, 0).
      \end{array}\right .
  \end{equation}

It is readily verified that $\Psi_\beta(s)$   possesses constant jumps
  $$
\left (\Psi_\beta\right )_+(s)=\left (\Psi_\beta\right )_-(s)
  \left\{
  \begin{array}{ll}
  \left(
    \begin{array}{cc}
      0 & 1 \\
      -1 & 0 \\
    \end{array}
  \right), &s\in (-\infty, 0),\\[.4cm]
  \left(
  \begin{array}{cc}
  -1 & 0 \\
  - 1  & -1 \\
  \end{array}
\right),
   & s\in (0, \infty) ,
  \\[.4cm]
 \left(
    \begin{array}{cc}
      1 & -1\\
     0 & 1 \\
    \end{array}
  \right), &s\in  i(0,\infty), \\ [.4cm]
   \left(
    \begin{array}{cc}
      1 & 1 \\
       0 & 1 \\
    \end{array}
  \right), &s\in  -i(0,\infty),
     \end{array}\right .
  $$
where orientation is taken from left to right on the real line, and down to up along the imaginary axis.

Now we introduce   a conformal mapping
\begin{equation}\label{beta-conformal-mapping}
  \lambda_\beta(z)= \left ( \frac 3 2\right )^{2/3} n^{1/3} \left (-\phi_\beta (z)\right )^{2/3}
\end{equation} from an $O(1)$ neighborhood of $\tau_\beta=\frac {z-\beta}{\beta-1}$ of the origin to  an $O(1)$ neighborhood of $\lambda_\beta=0$.
To verify the fact, we may need a convergent series expansion of the form
\begin{equation}\label{phi-beta-expansion}
  n\phi_\beta (z)\sim -c_\beta(n)\tau_\beta^{3/2},~~\arg\tau_\beta\in (-\pi, \pi),
\end{equation}where $\tau_\beta=\frac {z-\beta}{\beta-1}$,  and $c_\beta(n) \sim \frac {2\sqrt{2b}} 3 \sqrt n$ for large $n$.

The parametrix can be represented as
\begin{equation}\label{P-beta-representation}
  P^{(\beta)}(z)=E_\beta (z) \Psi_\beta \left (n^{1/3}\lambda_\beta(z)\right )e^{n\phi_\beta(z)\sigma_3}e^{\pm \frac \pi 2 i\sigma_3}, ~~\pm\Im z>0,~~z\in U_r(\beta),
\end{equation}where $E_\beta(z)$ is an analytic function in $z\in U_r(\beta)$. A straightforward verification shows  that, for $s=n^{1/3}\lambda_\beta(z)$, $\zeta=\frac 2 3 s^{\frac 3 2}=-n\phi_\beta(z)$, and $\zeta=-n\phi_\beta(z)\sim \sqrt n \tau_\beta^{3/2}\xi_0(\tau_\beta)$ for mild $\tau_\beta$ and large $n$,
 \begin{equation}\label{Psi-beta-expansion}
    \Psi_\beta \left (s\right )e^{n\phi_\beta(z)\sigma_3}e^{\pm \frac \pi 2 i\sigma_3}\sim
    \frac {s^{-\frac 1 4 \sigma_3}}{2\sqrt \pi}  \left\{ \sum_{k=0}^\infty \frac {\Psi_{k}}{\zeta^k}\right\} M_\beta,~~\arg z\in (-\pi, \pi);
\end{equation}
 see \cite[(10.4.7), (10.4.59), (10.4.61)]{as}, where $M_\beta=\left(
       \begin{array}{cc}
         1 & i \\
         1 & -i \\
       \end{array}
     \right)$, and the coefficients $\Psi_{k}$ are given in \eqref{Psi-expansion}.

Denoting  $f(\tau)=\frac {c_1 \sqrt{2/b}}{\tau^2\xi_0(\tau)}$, and $f_s(\tau)=\frac 1 {12 b} \frac 1 \tau+ \frac 5 {48b} \frac 1 {\tau^2}$ being  the  singular part of it,
we rewrite $N(z)$ as
$$N(z)= \left ( I+O\left (\frac 1 n\right )\right ) M_\beta^{-1} \left(
                       \begin{array}{cc}
                         1 & f_s(\tau_\beta) \\
                         0 & 1 \\
                       \end{array}
                     \right)M_\beta   M_\beta^{-1} \varrho(z)^{-\sigma_3} M_\beta.$$
   Thus,   we may choose
\begin{equation}\label{E-beta}
  E_\beta(z)= 2\sqrt\pi M_\beta^{-1} \left(
                       \begin{array}{cc}
                         1 & -f_r(\tau_\beta) \\
                         0 & 1 \\
                       \end{array}
                     \right) \left ( \frac {s^{\frac 1 4}}{ \varrho(z)}\right )^{\sigma_3}   ,~~z\in U_r(\beta),
\end{equation}where again $f_r(\tau)=f(\tau)- f_s(\tau)$ is the regular part of $f(\tau)$. With the analytic factor so chosen, it is readily verified that the matching condition is satisfied.

For later use, we write down the well-known formulas
\begin{equation}\label{Airys}
  \omega\Ai(\omega s)=-\frac 1 2 \left (\Ai(s)-i\Bi(s)\right )~~\mbox{and}~~\omega^2\Ai(\omega^2 s)=-\frac 1 2 \left (\Ai(s)+i\Bi(s)\right )~~\mbox{for}~s\in \mathbb{C}.
\end{equation}

\subsection{The final transformation $S\to R$}
\begin{figure}[t]
\begin{center}
\includegraphics[width=12cm,bb=29 315 545 474]{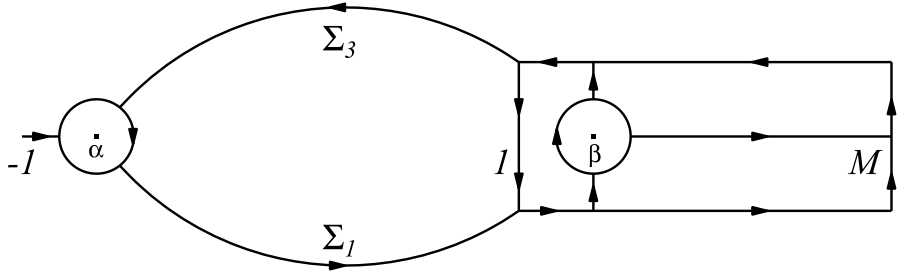}
 \caption{The contour $\Sigma_R$.}
 \label{contour-for-R}
\end{center}
\end{figure}

We bring in the final transformation by defining
\begin{equation}\label{R-define}
R(z)=\left\{ \begin{array}{ll}
                S(z)N^{-1}(z), & z\in \mathbb{C}\backslash \left ( U_r(\beta)\cup U_r(\alpha)\cup \Sigma_S \right ), \\
               S(z) (P^{(\beta)})^{-1}(z), & z\in   U_r(\beta)\backslash \Sigma_S ,  \\
               S(z)  (P^{(\alpha)})^{-1}(z), & z\in   U_r(\alpha)\backslash
               \Sigma_S .
             \end{array}\right .
\end{equation}
We note that the jumps along contours emanating from $z=1$ is of no significance; cf. \eqref{S-jump}, \eqref{N-jump}, and the discussion in Section \ref{phi-alpha-phi-beta}. On the remaining contours illustrated  in Figure \ref{contour-for-R},  we have the jump $J_R(z)=O(1/  n)$ uniformly for large $n$. Hence we conclude that $R(z)=I+O(1/  n)$.

\section{Proof of Theorem \ref{Main-theorem}}
\setcounter{equation} {0}  \label{proof-of-theorem}

We prove the theorem case by case, tracing back to the transformations  $Y\to U\to T \to S \to R$. The regions are illustrated in Figure \ref{region-asymptotic}.

From \eqref{Y-def} it is known that
$
\pi_n(z)=Y_{11}(z)$, the $(1,1)$ entry of $Y(z)$.
Also, it is readily seen from \eqref{pn} and the initial condition that the leading   coefficient of the orthonormal polynomial $p_n(z)$ is   $\gamma_n=2^n$, namely
\begin{equation}\label{proof-pn}
p_n(z)=\gamma_n\pi_n(z)=2^n\pi_n(z).
\end{equation}

For $z\in A_r$, we have
 \begin{equation*}
   Y(z)=e^{\frac 12nl\sigma_3}R(z)N(z)e^{n(g(z)-\frac 12l)\sigma_3}\{d_E(z)\chi(z)\}^{-\sigma_3}.
\end{equation*}
 Accordingly,
\begin{equation*}
\pi_n(z)=e^{ng(z)}\{d_E(z)\chi(z)\}^{-1}(R_{11}N_{11}+R_{12}N_{21}).
\end{equation*}
In view  of the fact that  $R(z)=I+O(1/n)$ as $n\rightarrow\infty$,
from \eqref{d-E}, \eqref{chi}, \eqref{phi-beta-g-relation}, \eqref{N-expression} and \eqref{proof-pn},
 we  obtain \eqref{thm-pi-ar}.

For $z\in B_r$, it is easily  seen  that
\begin{equation}\label{Y-br}
Y(z)=e^{\frac 12nl\sigma_3}R(z)N(z)
\begin{pmatrix}
1 & 0\\
\phi_0(z)e^{2n\phi_\alpha(z)} & 1
\end{pmatrix}
e^{n(g(z)-\frac 12l)\sigma_3}\{d_E(z)\chi(z)\}^{-\sigma_3}.
\end{equation}
From \eqref{Y-br}
we have
\begin{align*}
\pi_n(z)=&e^{ng(z)}\{d_E(z)\chi(z)\}^{-1}[(R_{11}N_{11}+R_{12}N_{21})
+(R_{11}N_{12}+R_{12}N_{22})\phi_0(z)e^{2n\phi_\alpha(z)})]\nonumber\\
=&e^{ng(z)+n\phi_\alpha}\phi_0^\frac{1}{2}(z)\{d_E(z)\chi(z)\}^{-1}
[R_{11}(N_{11}e^{-n\phi_\alpha}\phi_0^{-\frac{1}{2}}+
N_{12}e^{n\phi_\alpha}\phi_0^{\frac{1}{2}})\nonumber\\
& +R_{12}(N_{21}e^{-n\phi_\alpha}\phi_0^{-\frac{1}{2}}
+N_{22}e^{n\phi_\alpha}\phi_0^{\frac{1}{2}})].
\end{align*}
%It follows from \eqref{N-expression} that
%\begin{align}\label{N11-N12-phi}
%&N_{11}e^{-n\phi_\alpha}\phi_0^{-\frac{1}{2}}+
%N_{12}e^{n\phi_\alpha}\phi_0^{\frac{1}{2}}\nonumber\\
%&=D^{-1}(\infty)\left[\left(1+f_s(\tau_\beta)\right)\varrho e^{-\frac{1}{4}\pi i}\cos \Theta_B+\left(1-f_s(\tau_\alpha)-f_s(\tau_\beta)f_s(\tau_\alpha)\right)
%\varrho^{-1} e^{\frac{1}{4}\pi i}\sin \Theta_B\right]
%\end{align}
%and
%\begin{align}\label{N21-N22-phi}
%&N_{21}e^{-n\phi_\alpha}\phi_0^{-\frac{1}{2}}
%+N_{22}e^{n\phi_\alpha}\phi_0^{\frac{1}{2}}\nonumber\\
%&=D(\infty)\left[\left(1-f_s(\tau_\beta)\right)\varrho e^{\frac{1}{4}\pi i}\cos \Theta_B+\left(1+f_s(\tau_\alpha)-f_s(\tau_\beta)f_s(\tau_\alpha)\right)
%\varrho^{-1} e^{-\frac{1}{4}\pi i}\sin \Theta_B\right],
%\end{align}
%where $\Theta_B=in\phi_\alpha+\frac{1}{2} i\ln \phi_\alpha(z)-i\ln D(z)+\frac{1}{4}\pi$.
Then \eqref{thm-pi-br} follows from   \eqref{phi-0},  \eqref{d-E}, \eqref{chi}, \eqref{phi-alpha-g-relation},     \eqref{N-expression},  \eqref{proof-pn},  and the fact that $f_s(\tau_\beta)$ and $f_s(\tau_\alpha)$ are bounded.

For $z\in C_r$, we have
\begin{align}\label{Y-cr}
Y(z)=e^{\frac 12nl\sigma_3}R(z)N(z)
\begin{pmatrix}
1 & 0\\
-e^{2n\phi_\beta(z)-2\pi i\gamma(z)} & 1
\end{pmatrix}
e^{n(g(z)-\frac 12l)\sigma_3}\{d_I(z)\chi(z)\}^{-\sigma_3}(D_+^u)^{-1}(z).
\end{align}
Recalling  the relation between $g(z)$ and $\phi_\beta(z)$, the definitions of $r(z)$  and $D_+^u(z)$, we obtain
\begin{align}\label{proof-pi-cr}
\pi_n(z)=&e^{\frac{1}{2}nl-\frac 12\pi i}{(-\gamma'(z)w_d(z))}^{-\frac 12}
[(R_{11}N_{11}+R_{12}N_{21})e^{-n\phi_\beta(z)+\pi i\gamma(z)}\nonumber\\
&~~~-(R_{11}N_{12}+R_{12}N_{22})e^{n\phi_\beta(z)-\pi i\gamma(z)}].
\end{align}
Therefore \eqref{thm-pi-cr} follows from \eqref{N-expression}, \eqref{proof-pn},   the fact that
$-\gamma'(z)w_d(z)=4b\,e^{-2b}(1+O(1/n))$,  and that both $f_s(\tau_\alpha)$ and $f_s(\tau_\beta)$ are bounded.
%we obtain
%\begin{align}\label{proof-pi-cr}
%\pi_n(z)&=e^{\frac{1}{2}nl}{(-\gamma'(z)w_d(z))}^{-\frac 12}\left[
%R_{11}(N_{11}e^{-n\phi_\beta(z)+\pi i\gamma(z)+\frac{1}{2}\pi i}+
%N_{12}e^{n\phi_\beta(z)-\pi i\gamma(z)-\frac{1}{2}\pi i})\right.\nonumber\\
%&~~~+\left.R_{12}(N_{21}e^{-n\phi_\beta(z)+\pi i\gamma(z)+\frac{1}{2}\pi i}+
%N_{22}e^{n\phi_\beta(z)-\pi i\gamma(z)-\frac{1}{2}\pi i})\right],
%\end{align}
%At the moment, $-\gamma'(z)w_d(z)=4b\,e^{-2b}(1+O(1/n))$.

%Indeed, Note that
%\begin{align}\label{Y-cr}
%Y(z)=e^{\frac 12nl\sigma_3}R(z)N(z)
%\begin{pmatrix}
%1 & 0\\
%-e^{2n\phi_\beta(z)-2\pi i\gamma(z)} & 1
%\end{pmatrix}
%e^{n(g(z)-\frac 12l)\sigma_3}\{d_I(z)\chi(z)\}^{-\sigma_3}(D_+^u)^{-1}(z).
%\end{align}
%Hence, similar to \eqref{N11-N12-phi} and \eqref{N21-N22-phi}, we can derive \eqref{thm-pi-cr} from \eqref{D-upper}, \eqref{phi-beta-g-relation}, \eqref{r-def}, \eqref{d-I}, \eqref{chi}, \eqref{Y-cr} and the fact that $\frac{1-f_s(\tau_\beta)}{1+f_s(\tau_\beta)}$ the real part of whose denominator is positive is bounded.

For $z\in D_{1,r}\cup  D_{2, r}$,
the series of transformations again applies, for $z\in D_{1,r}$, we have
\begin{equation}\label{Y-D1}
Y(z)=e^{\frac 12nl\sigma_3}R(z)P^{(\beta)}(z)\begin{pmatrix}
1 & 0\\
-e^{2n\phi_\beta-2\pi i\gamma} & 1
\end{pmatrix}e^{n(g(z)-\frac 12 l)\sigma_3}\{d_I(z)\chi(z)\}^{-\sigma_3}(D_+^u)^{-1}(z).
\end{equation}
We can also rewrite  $\Psi_\beta(s)$   in \eqref{Psi-beta} as
$$
\Psi_\beta(s)=\begin{pmatrix}
\Ai(s) & \Bi(s)\\
\Ai'(s) & \Bi'(s)
\end{pmatrix}\begin{pmatrix}
-\frac 12 & -\frac 12\\
-\frac i2 & \frac i2
\end{pmatrix}.
$$
Combining the definitions of $d_I(z)$ in \eqref{d-I}, $\chi(z)$ in \eqref{chi}, $D_+^u(z)$ in \eqref{D-upper} and the relation between $g(z)$ and $\phi_\beta(z)$ yields
\begin{align}\label{Y-D1-cos-sin}
&\begin{pmatrix}
-\frac 12 & -\frac 12\\
-\frac i2 & \frac i2
\end{pmatrix}\begin{pmatrix}
1 & 0\\
e^{-2\pi i\gamma} & 1
\end{pmatrix}e^{\frac 12\pi i\sigma_3}e^{2n\phi_\beta\sigma_3}e^{n(g(z)-\frac 12 l)\sigma_3}\{d_I(z)\chi(z)\}^{-\sigma_3}(D_+^u)^{-1}(z)\nonumber\\
&=\begin{pmatrix}
-\cos(\pi\gamma) & *\\
\sin(\pi\gamma) & *
\end{pmatrix}(-\gamma'(z)w_d(z))^{-\frac 12\sigma_3}
\end{align}
Substituting \eqref{Y-D1-cos-sin} into \eqref{Y-D1} gives
 \begin{align}\label{pi-in-P-beta}
    \pi_n(z)&=\sqrt{\pi}\,e^{\frac 12nl}(-\gamma'(z)w_d(z))^{-\frac 12}\left\{R_{11}(z)\,[s^{\frac 14}\varrho^{-1} A_1(z)+s^{-\frac 14}\varrho \,(-f_r(\tau)+1) A_2(z)]\right.\nonumber\\
    &~~~\left.+R_{12}(z)\,[-is^{\frac 14}\varrho^{-1} A_1(z)+s^{-\frac 14}\varrho \,(if_r(\tau)+i) A_2(z)]\right\}.
    \end{align}
    Then \eqref{pi-n-beta} follows from \eqref{proof-pn} and the fact that
$-\gamma'(z)w_d(z)=4b\,e^{-2b}(1+O(1/n))$.
The case $z\in D_{2,r}$ can be treated similarly.
%Still we include some details. This time
%$$
%Y(z)=e^{\frac 12nl\sigma_3}R(z)P^{(\beta)}(z)e^{n(g(z)-\frac 12 l)\sigma_3}\{d_E(z)\chi(z)\}^{-\sigma_3}(D_+^l)^{-1}(z),
%$$
%and  $\Psi_\beta(s)$ can be written as
%$$
%\Psi_\beta(s)=\begin{pmatrix}
%\Ai(s) & \Bi(s)\\
%\Ai'(s) & \Bi'(s)
%\end{pmatrix}\begin{pmatrix}
%-\frac 12 & -1\\
%-\frac i2 & 0
%\end{pmatrix}.
%$$
%Also, it is readily verified that
%the functions $d_E(z)$, $\chi(z)$, $g(z)$ and $\phi_\beta(z)$ obey the relation
%$$
%e^{\frac 12\pi i\sigma_3}e^{2n\phi_\beta\sigma_3}e^{n(g(z)-\frac 12 l)\sigma_3}\{d_E(z)\mathcal{X}(z)\}^{-\sigma_3}=e^{\pi i\gamma\sigma_3}(-\gamma'(z)w_d(z))^{-\frac 12\sigma_3},
%$$and that
%$D_+^l(z)$ in \eqref{D-lower} fulfils
%$$
%(-\gamma'(z)w_d(z))^{-\frac 12\sigma_3}(D_+^l)^{-1}(z)=\begin{pmatrix}
%\frac{2i\sin(\pi\gamma)}{e^{i\pi \gamma}} & 0\\
%1 & \frac{e^{i\pi\gamma}}{2i\sin(\pi\gamma)}
%\end{pmatrix}(-\gamma'(z)w_d(z))^{-\frac 12\sigma_3},
%$$
%and that the following identity holds
%$$
%\begin{pmatrix}
%-\frac 12 & -1\\
%-\frac i2 & 0
%\end{pmatrix}e^{\pi i\gamma\sigma_3}\begin{pmatrix}
%\frac{2i\sin(\pi\gamma)}{e^{i\pi \gamma}} & 0\\
%1 & \frac{e^{i\pi\gamma}}{2i\sin(\pi\gamma)}
%\end{pmatrix}=\begin{pmatrix}
%-\cos(\pi\gamma) & -\frac{1}{2i\sin(\pi\gamma)}\\
%\sin(\pi\gamma) & 0
%\end{pmatrix}.
%$$
%Combining together all above now yields \eqref{pi-in-P-beta} for $z\in D_{2,r}$. Thus \eqref{pi-n-beta} follows accordingly.

For $z\in E_r$, it is readily seen that
\begin{align}\label{Y-er}
Y(z)=e^{\frac 12nl\sigma_3}R(z)N(z)e^{n(g(z)-\frac 12l)\sigma_3}\{d_I(z)\chi(z)\}^{-\sigma_3}(D_+^l)^{-1}(z).
\end{align}
Picking up the $(1,1)$-entry gives
\begin{align}\label{proof-pi-er}
\pi_n(z)=&e^{\frac{1}{2}nl-\frac 12\pi i}{(-\gamma'(z)w_d(z))}^{-\frac 12}
[2i\sin(\pi\gamma(z))e^{-n\phi_\beta(z)}(R_{11}N_{11}+R_{12}N_{21})\nonumber\\
&~~~-(R_{11}N_{12}+R_{12}N_{22})e^{n\phi_\beta(z)-\pi i\gamma(z)}].
\end{align}
Now a combination of the fact  that $D^{-1}(z)$ is bounded with \eqref{N-expression} and  \eqref{proof-pn}  justifies \eqref{thm-pi-er}.

For $z\in  F_{1, r}\cup  F_{2, r}$, first for $z\in F_{1, r}$, we have
\begin{equation}\label{}
  Y(z)=e^{\frac 12nl\sigma_3}R(z)P^{(\alpha)}(z)e^{n(g(z)-\frac 12 l)\sigma_3}\{d_E(z)\chi(z)\}^{-\sigma_3}.
\end{equation}
Then we have
\begin{equation}\label{pi-in-P-alpha}
\pi_n(z)=[R_{11}(z)(P^{(\alpha)})_{11}(z)+R_{12}(z)(P^{(\alpha)})_{21}(z)]e^{ng(z)}\{d_E(z)
\chi(z)\}^{-1}.
 \end{equation}
Therefore, substituting
  $P^{(\alpha)}(z)$ in \eqref{P-alpha-representation}, the relation    \eqref{phi-alpha-g-relation}, $\{d_E(z)\chi(z)\}^{-1}(\phi_0)^{1/2}=1+O(1/n)$ and  the asymptotic behavior of $R(z)$  into the representation,  and in view of  \eqref{proof-pn},  we   obtain \eqref{pi-n-alpha}. The case of $z\in F_{2, r}$ is similar to that of $z\in F_{1, r}$.

\section{Asymptotic quantities  and comparison with known results}
 \setcounter{equation} {0}  \label{comparison}
We evaluate several  asymptotic quantities and approximations  from the Riemann-Hilbert analysis and Theorem \ref{Main-theorem}, and compare them with the results of Wang and Zhao \cite{wang-zhao}  obtained earlier via integral methods.

 \subsection{Evaluation of the Lagrange multiplier $l$}
To calculate the constant $l=l(n)$ in \eqref{g-phase-condition}, and \eqref{phi-alpha-g-relation}-\eqref{phi-beta-g-relation}, we use the relations
\begin{gather}\label{g-jumps}
\begin{array}{ll}
   g_+(x)+g_-(x)-l +\frac 1 n  \ln r(x) =0,&  x\in (\alpha, \beta), \\
   g_+(x)-g_-(x)=2\pi i, & x\in(-\infty,\alpha),\\
   g_+(x)-g_-(x)=\frac  {2\pi i\,\gamma(x)} n,& x\in(\beta,\infty);
\end{array}  \end{gather}
cf. \eqref{g-phase-condition} and \eqref{g-in-phi}. It is also seen that
 $g(z)\sqrt{(z-\alpha)(z-\beta)}\rightarrow 0$ as $z\rightarrow \alpha$ and $z\to \beta$, and that
$g(z)=\ln z+O(z^{-1} \ln z)$ as $z\rightarrow\infty.$

An alternative representation is obtained by solving  the above scalar RHP, namely,
\begin{align}\label{g-integral}
g(z)&=\ln\frac{2z-\alpha-\beta+2\sqrt{(z-\alpha)(z-\beta)}}{\beta-\alpha}+
\frac{\sqrt{(z-\alpha)(z-\beta)}}{2n\pi}\int_\alpha^\beta\frac{\ln r(x)}{\sqrt{(x-\alpha)(\beta-x)}}\frac{dx}{x-z}\nonumber\\
&~~~+\frac {\sqrt{(z-\alpha)(z-\beta)}} n \int_\beta^\infty\frac{\gamma(x)}{ \sqrt{(x-\alpha)(x-\beta)}}
\frac{dx}{x-z}+\frac l2.
\end{align}
Let $z\rightarrow\infty$, we find that
\begin{equation}\label{l-representation}
\frac l2=\ln \frac{\beta-\alpha}{4}-\frac {\ln \pi}{2n}+\frac{1}{2n\pi}\int_\alpha^\beta\frac{ \ln (\pi  r(x))\;  dx }{\sqrt{(x-\alpha)(\beta-x)}}+\frac 1 n \int_\beta^\infty\frac{\gamma(x)\; dx}{\sqrt{(x-\alpha)(x-\beta)}}.
\end{equation}
To pick up the   contribution from the first integral,
we apply the Cauchy integral theorem to the function
 $\frac {-ib\ln\left (z+\sqrt{z^2-1}\right ) }{\sqrt{z^2-1}\sqrt{(z-\alpha)(z-\beta)}}$  on the cut-plane $\mathbb{C}\setminus (-\infty, \beta]$. As a result, we have
$$\int^1_\alpha \frac {\frac {2b\arccos x}{\sqrt{1-x^2}} dx }{\sqrt{(x-\alpha)(\beta-x)}}=
 \int^{-1}_{-\infty} \frac {2b\pi dx }{\sqrt{ x^2-1}\sqrt{( \alpha-x)(\beta-x)}} -\int_1^\beta \frac {2b\ln \left (x+\sqrt{x^2-1}\right ) dx }{\sqrt{ x^2-1}\sqrt{(x-\alpha)(\beta-x)}}.$$
     Recalling that $\ln (\pi r(x))=-\frac {2b\arccos x}{\sqrt{1-x^2}}$ for $x\in (\alpha, 1)$, and substituting the above equality of integrals into \eqref{l-representation}, we obtain
\begin{equation}\label{l-asymptotic}
nl=-2n\ln 2-\ln \pi+O(1/\sqrt n\; ).
\end{equation}Here, use has been made of the fact that $\alpha\sim -1+\frac b n$,     $\beta\sim  1+\frac b n$, the
the contribution of the integrals
  over the narrow intervals $[-1, \alpha]$ and $[1,\beta]$ are  negligible, and the leading terms of the integrals on $[\beta, +\infty)$ and $[-\infty, -1]$ can be evaluated and are  canceled with each other.  For example, since $\alpha\sim -1+\frac b n$, up to a factor $\left (1+O\left(\frac 1 n\right )\right )$, we have
$$\int^\infty_\beta \frac {dx}{\sqrt{x^2-1}\sqrt{(x-\alpha)(x-\beta)}}\sim \int^\infty_\beta \frac { dx}{(x+1) \sqrt{(x -1)(x-\beta)}}=\frac  {1 }
{\sqrt{2(\beta+1)}} \ln \frac {\sqrt{\frac {\beta+1} 2} +1}   {\sqrt{\frac {\beta+1} 2} -1}.$$ To deal with the second integral, a change of variable $t=\sqrt {\frac {x-\beta}{x-1}}$ may be used.

An alternative way to find $l$ is to use a certain Szeg\"{o} function, similar to \cite[(6.6)]{zhou-zhao}.

\subsection{Asymptotics of the leading  coefficients $\gamma_n$, a consistency check}
Write  $ Y(z)=\left(I+Y_1/z+\cdots\right)z^{n\sigma_3} $, with $Y_1$ independent of $z$.  Then
$
 \gamma_{n-1}^{2}=-\frac 1{2\pi i}(Y_1)_{21};
$ cf. \eqref{Y-def}, where $\gamma_{n-1}$ is the leading coefficient for the orthonormal polynomial  $p_{n-1}(z)$.
Tracing back the   series of transformations $Y\to U\to T\to S$ in the last section, we have
\begin{equation}\label{Y-in-S}
 Y(z)=e^{\frac 12nl\sigma_3}S(z)e^{n(g(z)-\frac 12l)\sigma_3}\{d_E(z)\chi(z)\}^{-\sigma_3}
\end{equation}
for $z$ in a neighborhood of infinity.

Elementary   calculation gives
\begin{equation}\label{d-E-chi}
    d_E(z)\chi(z)=1-b \frac{\ln (-z)}z+O\left (\frac 1 z\right )~~\mbox{as}~~z\to\infty;
\end{equation}cf. \eqref{d-E} and \eqref{chi}. Also, we have
\begin{equation}\label{g-at-infty}
  g(z)=\ln z-\frac bn\,\frac{\ln (-z)}{z}+O\left (\frac 1 z\right )~~\mbox{as}~~z\to\infty,
\end{equation}which follows from the fact that
$$G(z)=-\frac 1 {\pi i z}-\frac b { n\pi i} \frac {\ln (-z)}{z^2}+O\left( z^{-2}\right),$$ as can be seen from \eqref{G-solution}. In the previous formulas, the logarithms take principal branches, namely, $\arg t\in (-\pi, \pi)$ for $\ln t$.

We may set $S(z)=I+S_1/z+\cdots$. Then, substituting \eqref{d-E-chi} and \eqref{g-at-infty} into \eqref{Y-in-S} yields
$$
  \gamma_{n-1}^{2}=\frac i{2\pi }e^{-nl}(S_1)_{21}\sim \frac i{2\pi }e^{-nl}(N_1)_{21},
$$where $N(z)=I+N_1/z+\cdots$. Taking \eqref{N-expression} into account, it is readily seen that
$
(N_1)_{21}= -\frac i2+O(1/n)
$for large $n$.
Now we put together  the asymptotic approximation  \eqref{l-asymptotic} of the constant $l$, we obtain
\begin{equation}\label{gamma-n}
 \gamma_{n-1}\sim 2^{n-1}
\end{equation} for large $n$, which agrees with the fact that $\gamma_n= 2^n$, as  can be easily verified from  \eqref{pn} and the initial conditions attached.

\subsection{Comparison  with  \cite{wang-zhao}: The behavior at $z=\alpha$ and $z=\beta$ }

In Wang and Zhao \cite{wang-zhao}, using integral methods, uniform asymptotic approximations have been obtained in an interval of width $O(1/n)$ to which  $\alpha$ belongs.
More precisely, denoting  $z=-\cos\theta'=-\cos(t/\sqrt
n)$, the leading term of the asymptotic expansion
\begin{align}\label{uniform-asymptotic-at-alpha}
p_n(z)&=\frac{(-1)^n e^{\pi b/\theta'-b}}{\sqrt{2b}}\left(\frac{B^2(t)}{2b-t^2}\right)^{1/4}\left\{
\Ai(n^{1/3} B^2)n^{1/3}(1+O(n^{-1/2}))\right.\nonumber\\
&~~~\left.-\Ai'(n^{1/3}B^2)O(n^{-1/3})\right\},~~~n\rightarrow\infty
\end{align} holds   uniformly  for $\delta\leq t\leq M$, with $\delta$ and $M$ being
  positive constants, where $B^2(t)$ is an analytic function at $t=\sqrt{2b}$ with $B^2(t)=-2(2b)^{-1/6}(t-\sqrt{2b})+\cdots$; see \cite[(5.9)]{wang-zhao}.
    % \begin{equation}\label{Bt}
%B(t)=\left\{
% \begin{array}{ll}
 % \left [\frac{3}{2}\left (\frac{2b}{t}\arccos\left (\frac t {\sqrt{2b}}\right )
  %-\sqrt{2b-t^2}\right )\right ]^{1/3},   & 0<t<\sqrt{2b};\\
  %  i\left [\frac{3}{2}\left (-\frac{2b}{t}\arccosh \left (\frac{t}{\sqrt{2b}}\right )+\sqrt{t^2-2b}\right )
  % \right ]^{1/3}, &
  % t>\sqrt{2b}.
% \end{array}
%\right.
%\end{equation}

We are in a position to compare \eqref{pi-n-alpha} with \eqref{uniform-asymptotic-at-alpha}, starting by showing that
$B^2$ serves as a  conformal mapping at $z=\alpha$, just as $\lambda_\alpha$ does; cf. \eqref{alpha-conformal-mapping}.
Indeed, for $z-\alpha=(1+\alpha)\tau$, the parameters are connected as  $t^2-2b=2b\tau+O(1/n)$, that is, $t-\sqrt{2b}=O(\tau)+O(1/n)$. Several facts are readily  observed:
\begin{eqnarray*}
 & &  \lambda_\alpha(z)= -(2b)^{1/3}\tau(1+O(\tau))(1+O(1/\sqrt{n}))
 %\sim -2(2b)^{-1/6}(t-\sqrt{2b})
 , \\
&&
% n(g(z)+\phi_\alpha(z))=-b+b\pi/\theta'-n\ln 2+ n\pi i+O(1/\sqrt{n}),
\sqrt{2(z+1)}=\theta'+O(1/n),\\
 && \varrho(z)=2^{\frac 14}b^{-1/4}n^{\frac 14}(-\tau)^{-\frac 14}(1+O(1/n)).
\end{eqnarray*}
 Substituting the approximations into \eqref{pi-n-alpha}, we have
\begin{align*}\label{pn-lead-term}
p_n(z)&=\frac{(-1)^ne^{-b+\pi b/\theta'}}{\sqrt{2b}}\left\{n^{1/3}
(2b)^{-\frac 16} \Ai(n^{1/3}\lambda_\alpha(z))(1+O(\tau))(1+O(n^{-1/2}))\right.\nonumber\\
&~~~\left.-\Ai'(n^{1/3}\lambda_\alpha(z))O(n^{-1/3})\right\},
\end{align*}
which  agrees  with \eqref{uniform-asymptotic-at-alpha}. %Thus $x=\alpha$ serves as an turning point in a sense that $\pi_n(x)$ is exponential on the left, and oscillating on the right.
%\vskip .5cm

In Wang and Zhao \cite{wang-zhao}, a full uniform asymptotic expansion   at $z=\beta$  is derived. For $z=\cos(t/\sqrt n)$, the leading behavior  of the asymptotic expansion, in terms of the  variables in this paper,  reads
\begin{equation}\label{uniform-asymptotic-at-beta}p_n(z)
=\frac{e^b}{\sqrt{2b}}\left(-\frac{B^2(t)}{t^2+2b}\right)^{1/4}\widetilde A(n^{1/3}B^2) n^{1/3}(1+O(n^{-1/2}))- \widetilde A'(n^{1/3}B^2)O(n^{-1/3})
\end{equation} as $n\rightarrow\infty$,
 uniformly  for $t\in i[-M, -\delta]$ with positive  constants $M$ and $\delta$
satisfying $\delta<\sqrt{2b}<M$, where $\widetilde A(s)=-\cos(\pi \gamma(z))\Ai(s)+\sin(\pi \gamma(z))\Bi(s)$ and $B^2(t)$ are analytic function of $t$ at $t=-i\sqrt{2b}$ with $B^2(t)=2(2b)^{-1/6}(|t|-\sqrt{2b})+\cdots$.

For $z-\beta=(\beta-1)\tau$, we see that $|t|-\sqrt{2b}=O(\tau)+O(1/n)$. Quantities are related as
$$
\lambda_\beta(z)=(2b)^{1/3}\tau(1+O(\tau))(1+O(1/\sqrt{n})),~~ \varrho^{-1}(z)=
2^{\frac 14}b^{-\frac{1}{4}}n^{\frac 14}\tau^{-\frac 14}(1+O(1/n)).
$$
Hence, with the above arguments, we see that \eqref{pi-n-beta}  agrees  with \eqref{uniform-asymptotic-at-beta}, with $A_1$ and $A_2$ in \eqref{pi-n-beta} correspond to $\widetilde A$ and $\widetilde A'$ in \eqref{uniform-asymptotic-at-beta}.
\vskip .5cm

In Szeg\"{o}
\cite{Szego1975}, an  observation was made to  the Pollaczek polynomials that they show a singular behavior as compared with the classical polynomials.  It is worth noting that the sieved
Pollaczek polynomials also show  such a singular behavior  in some aspects,   such as the orders in $n$ of the behavior at the endpoints,  the Toeplitz minima, the behavior of $p_n(z)$ for $z$ not on the orthogonal support  or within  $(-1, 1)$, and the large-$n$ behavior of the  smallest zeros.

%For example, a combination of \eqref{pi-n-alpha} and  \eqref{l-asymptotic}
%$$ C(-1)^n e^{-n\ln 2} n^{\frac 1 3}.$$
%compared with the order $n^{\alpha+\frac 1 2}$ and $n^{\beta+\frac 1 2}$ at $x=\pm 1$.
%Toeplitz minima.
%If $x$ is not on the cut and nodes.
%$  (x+(x^2-1)^{\frac 1 2}  )^n$, the polynomial   $n^K    (x+(x^2-1)^{\frac 1 2}  )^n$, $K$ is a function of $x$.
%Discrepancy $x\in [-1, 1]$.
% Largest zeros.

\section*{Acknowledgements}
The work of Xiao-Bo Wu    was supported in part by the National Natural Science Foundation of   China under grant number 11201070, and the Educational Commission  of Guangdong Province  under grant number  Yq2013161.
 The work of Shuai-Xia Xu  was supported in part by the National
Natural Science Foundation of China under grant number
11201493, GuangDong Natural Science Foundation under grant number S2012040007824, and the Fundamental Research Funds for the Central Universities under grant number 13lgpy41.
 Yu-Qiu Zhao  was supported in part by the National 
Natural Science Foundation of China under grant number 
10871212.

\end{document}